\documentstyle{article}

\title{Interior estimates for solutions of Abreu's equation}
\author{S. K. Donaldson}
\date{\today}
\begin{document}
\maketitle

\newtheorem{thm}{Theorem}
\newtheorem{rmk}{Remark}
\newtheorem{prop}{Proposition}
\newtheorem{defn}{Definition}
\newtheorem{case}{Case}
\newtheorem{condn}{Condition}
\newtheorem{cor}{Corollary}
\newtheorem{lem}{Lemma}
\newcommand{\End}{{\rm End}}
\newcommand{\uA}{\underline{A}}
\newcommand{\cE}{{\cal E}}
\newcommand{\bR}{{\bf{R}}}
\newcommand{\bC}{{\bf{C}}}
\newcommand{\bZ}{{\bf{Z}}}
\newcommand{\utheta}{\underline{\theta}}
\newcommand{\Vol}{{\rm Vol}}
\newcommand{\Av}{{\rm Av}}
\newcommand{\db}{\overline{\partial}}
\newcommand{\ux}{\underline{x}}
\newcommand{\oz}{\overline{z}}
\newcommand{\ur}{\underline{r}}
\newcommand{\uSigma}{\underline{\Sigma}}
\newcommand{\uR}{\underline{R}}
\newcommand{\uv}{\underline{v}}
\newcommand{\un}{\underline{n}}
\newcommand{\uC}{\underline{C}}
\newcommand{\ukappa}{\underline{\kappa}}
\newcommand{\tkappa}{\tilde{\kappa}}
\newcommand{\tF}{\tilde{F}}
\newcommand{\dtheta}{\frac{\partial}{\partial \theta}}
\newcommand{\dH}{\frac{\partial}{\partial H}}
\newcommand{\dt}{\frac{\partial}{\partial t}}
\newcommand{\Euc}{{\rm Euc}}
\newcommand{\Riem}{{\rm Riem}}
\newcommand{\Ric}{{\rm Ric}}
\newcommand{\cS}{{\cal S}}
\newcommand{\partialone}{\partial_{x_{1}}}
\newcommand{\partialtwo}{\partial_{x_{2}}}
\newcommand{\partialxi}{\partial_{\xi}}
\newcommand{\partialeta}{\partial_{\eta}}
\newcommand{\ad}{{\rm ad}}
\newcommand{\diag}{{\rm diag}}
\newcommand{\normal}{{\rm norm}}
\newcommand{\Area}{{\rm Area}}

\section{Introduction}

In this paper we study a nonlinear, fourth order, partial differential
equation for a convex function $u$ on an open set $\Omega$ in $\bR^{n}$.
The equation can be written as
$$ S(u)= A$$

where $A$ is some given function and $S(u)$ denotes the expression
\begin{equation}  S(u)=- \sum_{i,j} \frac{\partial^{2} u^{ij}}{\partial x^{i}\partial x^{j}}
. \label{eq:Abreu}\end{equation} Here $(u^{ij})$ denotes the inverse of
the Hessian matrix $u_{ij}= \frac{\partial^{2} u}{\partial x^{i}\partial
x^{j}}$. We call this PDE {\it Abreu's equation} since the expression $S(u)$
 appears in \cite{kn:Abr}, in the study of the differential
geometry of toric varieties. In this paper we will be primarily interested
in the case when $A$ is a   constant. Solutions to Equation~\ref{eq:Abreu} then correspond to certain
 Kahler
metrics of constant scalar curvature, as we will recall in more detail
in Section 5 below. Our purpose  is to derive  {\it
a priori} estimates for solutions of Abreu's equation, which can be applied
to existence questions for such constant scalar curvature metrics, on the
lines of \cite{kn:don}.  However in the present paper we will 
 keep the differential
geometry in the background, and concentrate on the PDE aspects of the
 equation.
 
 Abreu's equation can be fitted, as a limiting case, into a class of equations
 considered by Trudinger and Wang \cite{kn:TW1},\cite{kn:TW2}. These authors study  the Euler-Lagrange
 equation of the functional   
   $$  {\cal J}_{\alpha}(u)= \int_{\Omega} ( \det (u_{ij}))^{\alpha} ,$$
  (or this functional plus lower order terms). The Legendre transform interchanges the equations with
    parameters
   $\alpha$ and $1-\alpha$. There are two exceptional cases, when $\alpha=0$
   or $1$, when the Euler-Lagrange equations are trivial. Abreu's equation,
   which is the Euler-Lagrange equation associated to the functional
   $$   \int_{\Omega} -\log \det (u_{ij}) + A u , $$ 
   is the natural limit of the Trudinger-Wang family when $\alpha\rightarrow
   0$. Indeed, the Trudinger-Wang equations (in the absence of lower order
   terms) can be written as
   $$      \sum_{ij} \frac{\partial^{2}}{\partial x^{i}\partial x^{j}}
    (\det(u_{ij})^{\alpha} u^{ij}) = 0. $$

    We will study Abreu's equation augmented by some specific boundary
    conditions. These depend on a measure $\sigma$ on the boundary of $\Omega$.
    We will consider two cases
    \begin{case}
    $\Omega$ is the interior of a bounded polytope, defined by a finite number
    of linear inequalities with $n$ codimension-$1$ faces of $\partial
    \Omega$ meeting at each
    vertex. The measure $\sigma$ is a constant muptiple
    of the Lebesgue measure on each codimension-$1$ face of $\partial
    \Omega$.
    \end{case}
    \begin{case}
    $\Omega$ is a bounded domain with strictly convex smooth boundary and
    $\sigma$ is a smooth positive measure on $\partial \Omega$.
    \end{case}
  
  The first case is the one which is relevant to toric varieties and is
  our main concern. We include the second case because it seems to lead
  to a natural
  PDE problem. In either case we define a class of convex functions
  $\cS_{\Omega,\sigma}$ satisfying boundary conditions depending on $\sigma$.
  The detailed definitions are given in Section 2.2 below, but roughly
  speaking we require that a function $u$ in $\cS_{\Omega,\sigma}$ behaves like
  $ \sigma^{-1} d\log d$ where $d$ is the distance to the boundary and
  $\sigma$ is regarded as a function, i.e. a multiple of the area measure
  on $\partial \Omega$. (Note here that the whole theory is  affine-invariant,
  and does not depend on the choice of a Euclidean metric on $\bR^{n}$,
  but for simplicity we will sometimes, as just above, express things in
  terms of the metric structure, although this is not playing any real
  role.) In this paper we study solutions $u$ of Abreu's equation which
  lie in $\cS_{\Omega,\sigma}$. As explained in \cite{kn:don} and in Section
  2.2 below, these arise when one considers the problem of minimising the
  functional
    \begin{equation} {\cal F}(u) = - \int_{\Omega} \log \det u_{ij} + A u - \int_{\partial
    \Omega} u d\sigma,  \label{eq:functional}\end{equation}
    over the set of smooth convex functions on $\Omega$ with $L^{1}$ boundary
    values.
    It is likely (although this  requires proof) that any extrema
    of this functional lies in our space  $\cS_{\Omega,\sigma}$.

    As explained in \cite{kn:don}, simple examples show that for some data $\Omega, A, d\sigma$ there
    are no solutions of Abreu's equation with the given boundary behaviour.
    By the same token, one needs further assumptions before {\it a priori}
    estimates of a solution can be obtained. The condition which we expect
    to be appropriate involves the linear part of the functional ${\cal
    F}$;
    \begin{equation} {\cal L}(f) = \int_{\partial \Omega} f d\sigma - \int_{\Omega} f
    A d\mu. \label{eq:linfunct}\end{equation}
    Fix a base point in the interior of $\Omega$ and call a function normalised
    if it vanishes, along with its first partial derivatives, at this base
    point. We consider the following condition on the data $(\Omega, A, d\sigma)$:
    \begin{condn}
    The functional ${\cal L}$ vanishes on affine-linear functions and there
    is some $\lambda>0$ such that
    $$   {\cal L}(f) \geq \lambda^{-1} \int_{\partial \Omega} f d\sigma, $$
    for all convex normalised functions $f$ on $\overline{\Omega}$. \end{condn}
    
    Of course when this condition holds we can fix $\lambda=\lambda(\Omega,
    \sigma, A)$ by taking
    the best possible constant.
    It is shown in \cite{kn:don}, at least in the case when the dimension $n$ is $2$
    and $\Omega$ is a polygon, that this is a necessary condition for the
    existence of a solution. Our goal in this paper
     is to derive {\it a priori} interior estimates for
    solutions assuming this condition. To streamline the statement we introduce
    the following terminology. We say a function $C(\Omega,\sigma,A, \lambda,
    d)$ where $\lambda$ and $d$ are positive real variables is {\it tame}
    if it is continuous with respect to the natural topology on the space
    of variables $(\Omega, \sigma, A, \lambda,d)$ (We use the 
    $C^{\infty}_{{\rm loc}}\cap
    L^{\infty}$ topology on $A$; in Case 1 
    the space is divided into components labelled by the combinatorics
    of the  faces
    and in Case 2 we use the $C^{\infty}$ topology on $\Omega, \sigma$).
    We use the same terminology for functions that depend on some subset
    of the variables. 
      Our main result is 
    \begin{thm}
    Suppose the dimension $n$ is $2$. There are tame functions
    $K, C_{p}$ for $p=0,1,\dots$ such that if $A$ is a smooth bounded function
    in $\Omega$, and  $(\Omega, \sigma, A)$ satisfies Condition 1, then any
    normalised solution $u$ in $\cS_{\Omega,\sigma}$ of
 Abreu's equation 
    satisfies  $$  K^{-1}\leq   (u_{ij})\leq K$$
   and $$ \vert \nabla^{p} u \vert\leq C_{p} $$  where the argument $d$ is the distance
   to the boundary of $\Omega$ and the argument $\lambda$ is
   $\lambda(\Omega,\sigma,A)$ \end{thm}
   
   Of course we need not take the definition of \lq\lq tame''functions
   $C(\Omega,\sigma,A,\lambda,d)$  too seriously. With a little labour
   we could make all of  our estimates completely explicit. The definition
   is tailored to the continuity method: if we have a continuous 
   $1$-parameter family of such problems defined by data $(\Omega_{t},\sigma_{t}, A_{t})$
   with solutions for $t<t_{0}$ then provided $\lambda(\Omega_{t}, \sigma_{t},
   A_{t})$ stays bounded the solutions cannot blow up in the interior as
   $t\rightarrow t_{0}$.

   We present a variety of different arguments to establish these interior
   estimates. We hope this variety is justified by the desire to extend the
   results, in the future, in various directions: to  the behaviour near
   the boundary and to higher dimensions. In Section 2 we bring together
   some more elementary preliminaries and in particular show that
    Condition 1 gives $C^{1}$ bounds on the solution. After this, the 
     crucial intermediate goal is to obtain
   upper and lower bounds on the determinant of the Hessian $(u_{ij})$.
   In Section 3 we consider the case when the dimension $n$ is $2$ and
    $A$ is
    a constant. We find a special argument in this case using a
     property of solutions of
  general  elliptic equations in $2$ dimensions.  In Section 4
   we find other arguments, using the maximum
   principle, which apply in any dimension. We get a lower bound on the
   determinant in all cases and, using the technique of Trudinger and Wang,
   an upper bound involving also a 
   \lq\lq modulus of convexity''. As we explain in (5.1), in dimension
   $2$, this modulus of convexity is controlled by the lower bound on the
   determinant, using an old result of Heinz.
   Once we have established upper and lower
     bounds on the determinant we can appeal to
      sophisticated analysis of Caffarelli and Guti\'errez to complete the
      proof of Theorem 1. This is explained in (5.1). In the case when $A$ is a constant and $\Omega$
      is a polygon in $\bR^{2}$, which is our main interest, we give an
      alternative proof in the remainder of Section 5. This avoids the deep analysis of
      Caffarelli and Guti\'errez (and perhaps gives more explicit estimates)
      but uses their geometric results about the
      \lq\lq sections'' of a convex function in an essential way. The other
      ingredients are $L^{2}$ arguments, which make contact with Kahler
      geometry, a variant of Pogorelov's Lemma and standard linear theory
      (Moser iteration). We also explain that, in this case, one can avoid
      using Heinz's result, by combing results from Sections 3 and 4.
      In sum, in the case when $A$ is constant and $\Omega\subset \bR^{2}$ is
      a polygon we get a self-contained proof of Theorem 1, assuming only Chapter
      3 of \cite{kn:Gut}, and material from the textbook \cite{kn:GT}.

  \section{Preliminaries}
  \subsection{Miscellaneous formulae}  
     
This subsection consists of entirely elementary material. We will give
a number of useful equivalent forms of Abreu's equation, which are obtained
by straightfoward manipulation. Throughout we use traditional tensor calculus
notation, with summation convention.

    Begin by considering any convex function $u$ defined on an open set
    in $\bR^{n}$. We write $u_{ij}$ for the Hessian, $u^{ij}$ for its inverse
    and $U^{ij}$ for the matrix of cofactors i.e. $U^{ij} = \det(u_{ij})
    u^{ij}$.  We can associate to $u$ two second order, elliptic,
     linear differential
    operators
    \begin{equation}    P(f) = \left( u^{ij} f_{i}\right)_{j}, \label{eq:Pdef}
    \end{equation}
    \begin{equation}   Q(f) = \left(U^{ij} f_{i}\right)_{j}.
     \label{eq:Qdef}\end{equation}
     
    Let $L= \log \det(u_{ij})$. Then we have identities
    \begin{equation}    L_{i} = u^{ab} u_{abi}, \label{eq:id1} \end{equation}
    \begin{equation}    u^{jk}_{i}= - u^{ja} u_{ab i} u^{bk}. \label{eq:id2}
    \end{equation}
    Thus
    \begin{equation}     u^{ij} L_{j}= - u^{ij}_{j}. \label{eq:id3}\end{equation}
        This means, first,  that 
    $$  U^{ij}_{j}= \left( e^{-L} u^{ij}\right)_{j}=
     e^{-L}\left( u^{ij}_{j}
    - u^{ij} L_{j} \right)
=0. $$
Hence the operator $Q$ can be written as
\begin{equation}   Q(f) = U^{ij} f_{ij}. \label{eq:Qdef2} \end{equation}

Now the first form of  Abreu's equation, as in the Introduction is
$$  u^{ij}_{ij}= -A. $$
We define a vector field $v= (v^{j})$ by
\begin{equation}   v^{j}= -u^{ij}_{ i} \label{eq:vdef}\end{equation} 
So our second form of Abreu's equation is
\begin{equation}    v^{j}_{j} = A. \label{eq:Ab2}\end{equation}
(Of course, the left hand side of this expression is the ordinary divergence
of the vector field $v$.)
On the other hand, by Equation~\ref{eq:id3}, the vector field can also be expressed as
$$     v^{j}= u^{ij}L_{j}, $$
so we have our third form of Abreu's equation
\begin{equation}   \left(u^{ij} L_{i}\right)_{j} = A,
\label{eq:Ab3} \end{equation}
that is $$ P(L)=A.$$
Expanding out the derivative we have
$$   \left(u^{ij} L_{i}\right)_{j}= u^{ij} L_{,ij} + u^{ij}_{j} L_{i}=
   u^{ij}\left(L_{ij} - L_{i} L_{j}\right), $$
   so we get our fourth form of the equation
   \begin{equation}   u^{ij}( L_{ij} - L_{i}L_{j}) = A. \label{eq:Ab4}
   \end{equation}
   Finally, if we write 
   $$  F= \det(u_{ij})^{-1} = e^{-L}, $$
   then $$ u^{ij} F_{ij} = u^{ij}(L_{ij}-L_{i}L_{j}) e^{-L}= -Ae^{-L}, $$
   so we get our fifth form of the equation
   \begin{equation}  Q(F) = -A .\label{eq:Ab5}\end{equation}    
   
   \subsection{The boundary conditions.}
   
   We begin by giving a precise definition of the set $\cS_{\Omega,\sigma}$ of
   convex functions on $\Omega$. The definitions are different in
    the two cases. We start with Case 1, when $\Omega$ is a polytope.
    For any  point $P$ of $\partial \Omega$  we can choose affine co-ordinates
    $x_{i}$ on $\bR^{n}$ such that $P$ has co-ordinates $x_{i}=0, i=1,\dots
    n$ and a neighbourhood of $p$ in $\overline{\Omega}$ is defined by
     $p$ inequalities
    $$ x_{1}, x_{2}, \dots, x_{p} > 0. $$
    We can also choose the co-ordinates so that the normal derivative of
    $x_{i}$ on the face $x_{i}=0$ of the boundary is  $\sigma^{-1}$. 
    We call such co-ordinates {\it adapted to $\Omega$ at $P$}. 
    \begin{defn}
        In Case 1 the set $\cS_{\Omega,\sigma}$ consists of continuous convex
        functions $u$ on $\overline{\Omega}$ such that
        \begin{itemize}
       \item $u$ is smooth and  strictly convex in $\Omega$,
       \item The restriction of $u$ to each face of $\partial \Omega$ is
       smooth and strictly convex;
       \item In a neighbourhood of any point $P$ of $\partial \Omega$ the
       function $u$ has the form
       $$  u= \sum x_{i} \log x_{i} + f $$
       where $x_{i}$ are adapted co-ordinates, as above, and $f$ is smooth
       up to the boundary. 
    \end{itemize}
    \end{defn}
    
    (We say a smooth function is {\it strictly convex} if its Hessian is
    strictly positive.)
     
    Now turn to Case 2, when $\Omega$ has smooth boundary. If $P$ is a
    point of $\partial \Omega$ we can choose local co-ordinates 
    $\xi, \eta_{1}, \dots \eta_{n-1}$ near $P$ so that $\partial \Omega$
    is given by the equation $\xi=0$ and the normal derivative of $\xi$
    on the boundary is  $\sigma^{-1}$. Again, we call such co-ordinates
    adapted. Define functions $\alpha_{p}$ of a positive real variable
    by
    $$ \alpha_{1}(t) = - \log t \ ,\ \alpha_{2}(t)
    = t^{-1}\ , \alpha_{3}(t) = t^{-2}. $$
    
    \begin{defn}
    In Case 2 the set $\cS_{\Omega,\sigma}$ consists of continous convex functions
    $u$ on $\overline{\Omega}$ such that
    \begin{itemize}
    \item $u$ is   smooth and strictly convex on $\Omega$;
    \item In a neighbourhood of any point $P$ of $\partial \Omega$ there
    are adapted co-ordinates $(\xi,\underline{\eta_{i}})$ in which
    $$ u = \xi \log \xi + f$$
    where for $p\geq 1$, $p+q\leq 3$
    $$  \vert \nabla_{\xi}^{p} \nabla_{\underline{\eta}}^{q} f \vert=
     o(\alpha_{p}(\xi)).
    $$
    as $\xi\rightarrow 0$.
    \end{itemize}    
    \end{defn}
    
    Here the notation $\nabla^{p}_{xi} \nabla^{q}_{\underline{\eta}}$ means
    any partial derivative of order $p$ in the variable $\xi$ and total
    order $q$ in the $\eta_{i}$.

   Now, for small positive $\delta$, let $\Omega_{\delta}\subset \Omega$
   be the  set of points distance at least $\delta$ from $\partial
   \Omega$, i.e. a  \lq\lq parallel'' copy of $\partial \Omega$. 
   If $u\in \cS_{\Omega,\sigma}$ and 
   $\chi$ is a smooth function on $\Omega$, integration-by-parts over
   $\Omega_{\delta}$ gives the fundamental identity:
   \begin{equation} \int_{\Omega_{\delta}} u^{ij} \chi_{ij} =  \int_{\Omega_{\delta}}
    u^{ij}_{ij} \chi + \int_{\partial \Omega_{\delta}}-u^{ij}_{j} \chi
    + u^{ij} \chi_{j}. \label{eq:fundint}\end{equation}
    We can write the boundary terms as
    $$ \int_{\partial \Omega_{\delta}} \chi\ v_{\normal}  + \nabla_{X} \chi, $$
    where $v$ is the vector field introduced in (10) above, $v_{\normal}$ is
    its normal component and
    $X$ is the vector field, defined on a neighborhood of $\partial
    \Omega$ by 
    $$  X^{j} = u^{ij} \nu_{j}, $$
    $\nu_{j}$ being the unit normal to $\partial \Omega$ at the closest
    boundary point. (In Case 1, the vector field $X$ will be discontinuous
    near the \lq\lq corners'' of $\partial \Omega$ but this will
    not matter.)
    
    The main result we need is
    \begin{prop}
    In either case, if $u\in \cS_{\Omega,\sigma}$ then as $\delta\rightarrow 0$;
    $\vert X \vert = O(\delta) $
    and
    $ v_{\normal}$ converges uniformly to $ \sigma. $
    \end{prop}

   Here $\vert X\vert$ refers to the Euclidean length and we interpret
   $v_{\normal}$ as a function on $\partial \Omega$ in the obvious way, by
   taking the closest point of $\partial \Omega_{\delta}$.
   We assume Proposition 1 for the moment. Taking the limit as $\delta$
   tends to $0$ in Equation~\ref{eq:fundint}, we obtain
   \begin{cor}
   Let $u$ be in $\cS_{\Omega,\sigma}$ with $S(u)=-u^{ij}_{ij}\in 
   L^{\infty}(\Omega)$
   and let $\chi$ be a 
    continuous, convex function on $\overline{\Omega}$ smooth in the interior and
   with $\nabla \chi= o(d^{-1})$, where $d$ is the distance to  $\partial
   \Omega$. Then $u^{ij} \chi_{ij}$ is integrable in $\Omega$ and
   $$ \int_{\Omega} u^{ij} \chi_{ij} = \int_{\Omega} u^{ij}_{ij} \chi +
   \int_{\partial \Omega} \chi d\sigma.
   $$
   \end{cor}
   
   The main application we make of Corollary 1
   is the case when $\chi=u$. It is clear from the definitions that $\nabla
   u $ is $O(-\log d)$, hence $o(d^{-1})$, near the boundary, 
   and $u_{ij} u^{ij} = n$. So we have the identity
     $$ {\cal L} u = n $$
     where $A= -u^{ij}_{ij}$ and ${\cal L}$ is the linear functional
        defined in (3). Thus we obtain
     \begin{cor}
     Suppose $A\in L^{\infty}(\Omega)$ and that ${\cal L}$ satisfies
     Condition 1. Then if $u$ is a normalised function in $\cS_{\Omega,\sigma}$ which satisfies
     Abreu's equation $S(u)=A$ we have
     $$ \int_{\partial \Omega} u \ d\sigma \leq n\lambda . $$
     \end{cor}
     A simple  argument (\cite{kn:don}, Lemma 5.2.3) shows that the integral over the boundary,
     for a normalised convex function, controls the derivative in the interior
     and we have
     \begin{cor} Under the same hypotheses as Corollary 2, 
     $$ \vert \nabla u \vert \leq C d^{-n}, $$
     where $C$ depends tamely on $(\Omega, \sigma,\lambda)$.
     \end{cor}
     This first derivative bound is the seed which we wish to 
     develop in this paper to obtain bounds on higher derivatives. 
     
     We mention some other applications
     of Corollary 1. Here we will restrict for simplicity to the case when
     $\Omega$ is a polytope. 
     \begin{cor}
     In the case when $\Omega\subset \bR^{n}$ is a polytope, for any
     $u\in \cS_{\Omega,\sigma}$, $\vert\log
     \det u_{ij}\vert$ is integrable over $\Omega$ so
       for any $A\in L^{\infty}(\Omega)$
     the functional ${\cal F}_{A}$ is defined on $\cS_{\Omega,\sigma}$.
      The functional ${\cal F}_{A}$ is convex on $\cS_{\Omega, \sigma}$
      and the equation 
     $S(u)=A$ has at most one normalised solution $u$ in
      $\cS_{\Omega,\sigma}$. For such $u$,
     
     $$\int_{\Omega}  \vert \log \det(u_{ij}) \vert \leq C, $$
     where $C$ is a tame function of 
     $\Omega, \sigma, \Vert A\Vert_{L^{\infty}}, \lambda$. 
      \end{cor}
     
     The proof of this Corollary follows easily from the results
     in (\cite{kn:don}, subsections (3.3) and (5.1)). 
     The restriction to  Case 1  arises because in this
     case, as is
  will be clear from Proposition 5 below,  $S(u)$ is bounded for any
   $u\in \cS_{\Omega,\sigma}$.  In Case 2 we have not analysed the behaviour
   of $S(u)$ near the boundary, for general elements of
    $\cS_{\Omega, \sigma}$, 
     and we leave  
      this issue to be discussed elsewhere. 
      
     \subsection{Proof of Proposition 1.}
     
     In the course of the proof we will also establish some other
     properties of functions in $\cS_{\Omega,\sigma}$.
       Before beginning it is worth pointing that, in a
     sense the proofs of these boundary properties need not be taken too seriously at this stage of the
     development of the theory. At this stage we are in essence free to
     choose the definition of $\cS_{\Omega,\sigma}$ and we could impose any
     reasonable conditions we choose. For example we 
     we could take the conclusions of Proposition 1 as part
      of the {\it definition} of $\cS_{\Omega,\sigma}$.
     The discussion will only acquire an edge when one goes on to the {\it
     existence}
     theory. It is possible that the detailed definitions may need to be
     modified then.  At the present stage, the definitions we have concocted serve
     to indicate at least the  nature of the solutions we want to
     consider while being, we hope, sufficiently general to permit a sensible
     existence theory. 
     
     We begin with Case 1. This is not very different from the analysis
     of metrics on toric varities in \cite{kn:Abr}, \cite{kn:Guil}, but we include a discussion
     for completeness. We work in adapted co-ordinates around a boundary
     point, so
     $$ u = \sum_{i=1}^{p} x_{i} \log x_{i} + f$$
     where $f$ is smooth up to the boundary. By the second condition of
     Definition 1  we may
     choose the co-ordinates so that at $\ux=0$
     $$ \frac{\partial^{2} f}{\partial x_{i} \partial x_{j}} =\delta_{ij},
     $$
     for $p+1\leq i,j\leq n$. Thus 
     \begin{equation} \left( u_{ij}\right) = \diag(x_{1}^{-1}, \dots, x_{p}^{-1}, 1\dots,
     1) + \left( \psi_{ij}\right), \label{eq:Hessform}\end{equation}
     where $\psi_{ij}$ are smooth up to the boundary and $\psi$ vanishes
     at $\ux=0$ for $i,j\geq p+1$. 
     \begin{prop}
     \begin{itemize}
     \item $\det(u_{ij})= x_{1}\dots x_{p} \Delta$ where $\Delta$ is smooth up to
     the boundary and $\Delta(0)=1$.
     \item There are  functions $f_{i}$, $g_{ij}$, $h_{ij}$, all smooth
     up to the boundary and with $f_{i}(0)=1$, $h_{ij}(0)=0$, such that
     $$ \left( u_{ij}\right) = \diag( f_{1} x_{1}, \dots, f_{p} x_{p},
     f_{p+1}, \dots f_{n}) + \left( \sigma^{ij}\right), $$
     where $$  \sigma^{ij}= x_{i} x_{j} g_{ij}$$ for $1\leq i,j\leq p$;
     $$ \sigma^{ij} = x_{i} g_{ij} $$ for $1\leq i\leq p$, $j>p$;
     $$ \sigma^{ij}= h_{ij} $$ for $i,j>p$.
     \end{itemize}\end{prop}
     It is completely straightforward to verify Proposition 1, in Case
     1, given this Proposition 2. Notice that, according to  Proposition
     2, the matrix $u^{ij}$ is
     smooth up to the boundary.
     
     Let $\Lambda$ be the diagonal matrix
     $$ \Lambda= \diag( x_{1}^{1/2}, \dots x_{p}^{1/2}, 1\dots, 1). $$
     To prove the first item of Proposition 2 we write $H$ for the Hessian matrix $(u_{ij})$
     and consider the matrix $\Lambda^{2} H$. This has the form $1+ E$
     where the entries $E_{ij}$ of the matrix $E$ are smooth up to the boundary
     and $E_{ij}$ vanishes at $0$ except for the range $1\leq j\leq p,
     \ i>p$. Thus, at $0$, the matrix $E$ is strictly lower-triangular and
     so $\Delta=\det(1+E)$ is a smooth function taking the value  $1$ at $\ux=0$,
     and $$ \det(u_{ij})= \det\Lambda^{-2}\ \Delta = x_{1}^{-1}\dots x_{p}^{-1}
     \Delta. $$
     
     To prove the second item we consider the symmetric matrix $\Lambda H \Lambda$.
      We write
      $ \Lambda H \Lambda = 1+ \left( F_{ij}\right)$ and 
      Equation~\ref{eq:Hessform} yields
      $$  F_{ij}= A_{ij} x_{i}^{1/2} x_{j}^{1/2} $$
      for $1\leq i,j\leq p$;
      $$  F_{ij} = F_{ji}= B_{ij} x_{i}^{1/2}, $$
      for $1\leq i\leq p$ and $j>p$ where $A_{ij}, B_{ij}$ are smooth up
      to the boundary. In the remaining block, $p+1\leq i,j\leq n$, the
      $F_{ij}$ are smooth up to the boundary, vanishing at $0$. Let
      $\ad( \Lambda H \Lambda)$ be the matrix of co-factors, or 
      \lq\lq adjugate matrix'', so that
            $$ (\Lambda H \Lambda)^{-1} = \det(\Lambda H \Lambda)^{-1}
            \ad (\Lambda H\Lambda).
      $$
      We claim that $\ad(\Lambda H \Lambda)$ has the form $1+F'_{ij}$
       where $F'$ satisfies the same conditions as $F$ above; that is, 
      $$  F'_{ij}= A'_{ij} x_{i}^{1/2} x_{j}^{1/2} $$
      for $1\leq i,j\leq p$;
      $$  F'_{ij}=B'_{ij} x_{i}^{1/2}$$
      for $1\leq i\leq p$ and $j>p$ etc. The proof of this is  a straightforward
       matter
      of considering the various terms in the cofactor determinants which
      we leave to the reader. Given this, the second item of the Proposition       
       follows by writing
         $$\left( u^{ij}\right)= H^{-1} = \Delta^{-1}\  \Lambda\ \ad(\Lambda H \Lambda)
        \ \Lambda. $$

     We now turn to Case 2. For simplicity we will consider the case when
     $n=2$ (see the remarks at the beginning of this subsection). Moreover, to make the calculations easier, we will consider
     a special kind of adapted co-ordinate. Suppose the point $P$ is the
     origin and let the boundary of $\Omega$
     be represented by a graph $x_{1}= Q(x_{2})$ where $Q(0)=Q'(0)=0, Q''(0)<0$. Then
     we take the adapted co-ordinates
     $$ \eta=x_{2}\ ,\ \xi=\rho(\eta) ( Q(x_{2})-x_{1}), $$
     for a positive function $\rho$ equal to $\sigma^{-1}$, where
     $\sigma$ is regarded as a function of $x_{2}$ using the obvious
       parametrisation of the boundary.
       Without loss of generality
     we can  compute at  points
     on the $x_{2}$ axis,
     where $\eta=0$. We write $\partialone,\partialtwo$ for $\frac{\partial}{\partial
     x_{1}}, \frac{\partial}{\partial x_{2}}$ and $\partialxi,\partialeta$
     for
     $\frac{\partial}{\partial \xi}, \frac{\partial}{\partial \eta}$. Then,
     transforming between the co-ordinate systems $(x_{1}, x_{2})$ and
     $(\xi,\eta)$, we have
     $$   \partialone = \rho \partialxi\ ,\  \partialtwo = \partialeta+F\partialxi,
     $$
     where $F(\xi,\eta)$ has the form $F=\xi A(\eta)+B(\eta)$ and $B(0)=0,B'(0)<0$.
     Thus
     \begin{eqnarray}  \partialone^{2}&=& \rho^{2} \partialxi^{2},\\
      \partialone\partialtwo&=& \rho^{2} \partialxi\partialeta + \rho F
     \partialxi^{2} + \rho A \partialxi, \\
     \partialtwo^{2}&=& ( \partialeta^{2}+ 2F \partialxi\partialeta +
     F^{2} \partialxi^{2}) + (\xi A'+B'+A(\xi A+B)) \partialxi. \end{eqnarray}
     Applying this to a function $u=(\xi\log\xi-\xi) + f$ in $\cS_{\Omega,\sigma}$
     (where we have
     replaced $f$ by $f-\xi$ in Definition 2, which obviously makes no
     difference)
     we get
     $$ \left(u_{ij}\right) = 
     \left( \begin{array}{ll}\rho^{2} \xi^{-1} &  0\\ 0 & B'\log \xi \end{array}
      \right)  
   + \left( \begin{array}{ll}\alpha &\beta\\ \beta &\gamma\end{array}\right), $$
   where 
   $$ \alpha= \partialone^{2} f\ ,\  \beta=\rho A \log \xi +\partialone\partialtwo
   f\ ,\  \gamma = (\xi A'+B'+ A(\xi A+B)) \log \xi+ \partialtwo^{2}f. $$
   Then
    $$ \det(u_{ij})= (\rho^{2}B') \xi^{-1}\log \xi \ \Delta, $$
    where $$ \Delta= \left( 1+ \frac{\alpha \xi}{\rho^{2}}\right)\left(
    1+ \frac{\gamma}{B'\log \xi}\right) - \frac{\beta^{2}\xi}{\rho^{2}B'
    \log \xi}. $$
    Write
    $$ \lambda=\frac{\alpha \xi}{\rho^{2}}, \mu=\frac{\gamma}{B'\log \xi},
    \nu=\frac{\beta \xi}{ B' \rho^{2} \log \xi}. $$
    Then $$\Delta= (1+\lambda)(1+\mu) - \beta \nu , $$
    and the inverse matrix is
    $$ \left(u^{ij}\right) = \Delta^{-1}
     \left( \begin{array}{ll} \frac{\xi}{\rho^{2}} (1+\mu)
     &- \nu\\ -\nu & \frac{1}{B'\log \xi} (1+\lambda)\end{array} \right). $$
     Now first we evaluate at a point where $\eta=0$, so $B=0$. Then we have
     $$\alpha= \rho^{2} \partialxi\partialxi f, $$
     and by the definition of $\cS_{\Omega,\sigma}$ this is $o(\xi^{-1})$. So
     $\lambda $ is $o(1)$. We have
     $$ \beta= \rho^{2} \partialxi\partialeta f + \rho \xi A \partialxi\partialxi
     f + \rho A (\partialxi f + \log \xi), $$
     and applying the definition of $\cS_{\Omega,\sigma}$ to the various terms
     we see that $\beta=O(\vert \log \xi\vert)$. This means that
     $\nu$ is $O(\xi)$ and $\beta \nu$ is $ O(\xi \vert \log \xi\vert)$,
     which are both $o(1)$.
     Similarly, applying the definitions, we find that $\mu$ is $o(1)$.
      We see from this first that the vector field $Z$, which has
     components $u^{11}, u^{12}$ is $O(\xi)$, so verifying the first item
     of Proposition 1. We also see that
     \begin{equation} \det(u_{ij}) \sim \frac{\rho^{2}}{B'} \xi^{-1} \vert \log\xi\vert,
     \label{eq:asymdet2} \end{equation} since $\Delta \sim 1$.  To complete  the proof we need to differentiate
     again. We have
     $$ L = \log\det(u_{ij})= -\log \xi + \log(-\log \xi) + \log(-B') +
     \log \rho^{2} + \log \Delta, $$
     Thus, when $\eta=0$,
     \begin{eqnarray} \partialone L&=& \rho(\xi^{-1} + (\xi \log \xi)^{-1}) + \rho 
     (\frac{\partialxi\Delta}{\Delta}),\\
     \partialtwo L &=& \xi A \partialxi L + \frac{\partialeta \Delta}{\Delta}.
     \end{eqnarray}
     We claim that $\partialeta \Delta$ is $O(1)$ and $\partialxi\Delta$
     is
     $o(\xi^{-1})$. Given this, we have 
     $$ v^{1}= \Delta^{-1} ( \rho^{-1} (1+\mu) + o(\xi)), v^{2}=\Delta^{-1}
     (-\nu \rho \xi^{-1} + O(\vert\log \xi\vert^{-1}), $$
     and we see that the normal component $v^{1}$ converges to $\sigma=\rho^{-1}$
     as desired. To verify the claim we have to show that 
     $ \partialeta \lambda, \partialeta \mu, \partialeta (\beta \nu)$ are
     all $O(1)$ and $\partialxi \lambda,\partialxi \mu, \partialxi (\beta
     \nu)$ are all $o(\xi^{-1})$. This is just a matter of differentiating
     the formulae defining $\lambda, \mu,\nu$ and applying the definition
     of $\cS_{\Omega,\sigma}$ to each term: we omit the details.

   \section{The two dimensional case.}
   \subsection{The conjugate function}
   Throughout this Section 3 we suppose that $A$ is a constant and $u$ is
   a solution of Abreu's equation in $\Omega$.
   We begin  with second form of the equation
   $$  div(v)= v^{i}_{i} = A. $$
   The \lq\lq radial'' vector field $x^{i}$ on $\bR^{n}$ 
   has divergence $n$, so if we define a vector field $w$ by
   \begin{equation}   w^{i}= v^{i}- \frac{A}{n} x^{i}, \label{eq:wdef}\end{equation}
   then $w^{i}_{i}=0$: the vector field $w$ has divergence zero.
   Now define a function $h$ by
   \begin{equation}   h= u-u_{k} x^{k}. \label{eq:hdef}\end{equation}
   Then
   \begin{equation} h_{i} = u_{i} - u_{ki}x^{k} - u_{k}\delta^{k}_{i} = - u_{ki} x^{k},
   \label{eq:hderiv}\end{equation}
   so
   \begin{equation}   u^{ij}h_{i}= - x^{j}. \label{hderiv2}\end{equation}
   Thus, using Equation~\ref{eq:id3}, we have
   \begin{equation}  w^{i} = u^{ij} \tilde{L}_{j}, \label{eq:tLdef}\end{equation} 
   where $\tilde{L}= L+ \frac{A}{n} h$. 
   
   We now specialise (for the rest of this subsection) to the case when $n=2$.
   The special feature here is that divergence-free vector fields can be
   represented by Hamiltonians, so there is a function $H$ (unique
   up to a constant) with
   $$    w^{i}= \epsilon^{ij} H_{j}, $$
   where $\epsilon^{ij}$ is the skew-symmetric  tensor with $\epsilon^{12}=1$.
   \begin{lem}
   The function $H$ satisfies the equation $Q(H)=0$. \end{lem}
   To see this we write
   $$  \tilde{L}_{b}= - u_{bk} \epsilon^{kj} H_{j}. $$
   Then we have
   $$\tilde{L}_{ab}\epsilon^{ab}= 0$$ by the symmetry of second derivatives,
   whereas, differentiating Equation~\ref{eq:tLdef},
   $$\tilde{L}_{ab}\epsilon^{ab}= - \left(\epsilon^{ab} 
   u_{bk}\epsilon^{kj} H_{j}\right)_{a}= Q(H), $$
   since $U^{aj}= \epsilon^{ab} u_{bk} \epsilon^{kj}$.
   
   We  call $H$ the {\it conjugate function} since
   the relationship between the functions $H$ and $ \tilde{L}$, 
     is analogous to that between conjugate
   harmonic functions in two dimensions, except that they satisfy
   {\it different} linear elliptic equations: $P(\tilde{L})=0, Q(H)=0$.

   Next we use the boundary conditions, so we suppose that $u$
    is in $\cS_{\Omega,\sigma}$.
   We have shown in Proposition 1 
   that the normal component of the vector field $v$ on the boundary can
   be identified with given measure $\sigma$. Likewise, the vector field
   $ \frac{A}{2} x^{i}$ defines a (signed) measure $d\tau$ on the boundary.
   The condition that ${\cal L}$ vanishes on the constants implies that
     \begin{equation}
      \Vol(\partial \Omega, \sigma)= A \Vol(\Omega),\label{Aid}\end{equation}
   which in turn implies 
   $$  \int_{\partial \Omega} d\sigma - d\tau = 0. $$
   Thus there is a function $b$ on $\partial \Omega$, 
   unique up to a constant, with
   $db= \sigma-\tau$. 
   Observe that, in the case when $\Omega$ is a polygon, the function $b$
   is linear on each face of $\partial \Omega$. 
   For
   in this case $d\sigma$ is constant multiple of the Lebesgue measure
   on each face and the {\it normal} component of the
   radial field is also constant on each face. 
 
 \begin{lem} The function $H$ extends continously to $\overline{\Omega}$
 with boundary value $b$ (up to the addition of a constant). 
 \end{lem}
 First suppose that $H$ extends smoothly up to the boundary. The relation
 $  w^{j}= \epsilon^{ij} H_{i}$ asserts that the normal component of $w$
 is equal to the tangential derivative of $H$. But the normal component
 of $w$ is the derivative of $b$, so the derivative of $H-b$
 vanishes. For the general case we apply the same argument to a slightly
 smaller domain and take a limit, the details are straightforward.

 We can now apply a special result (\cite{kn:GT}, Lemma 12.6) for 
 solutions of elliptic equations
 $$ \sum_{i,j=1}^{2} a^{ij} \frac{\partial^{2} f}{\partial x^{i}\partial x^{j}}= 0
 $$ in two dimensions whose boundary values satisfy a \lq\lq three point''
 condition. Here we use Equation~\ref{eq:Qdef2} to see 
 that the equation $Q(H)=0$ has this form,
 with $a^{ij}=U^{ij}$. The three point condition requires that for any
 three points $X_{1}, X_{2}, X_{3}$ on $\partial \Omega$ the slope
 of the plane in $\bR^{3}$ containing $ (X_{i}, H(X_{i}))$ is bounded by
 some fixed $K$. This holds in our situation because the restriction of
  $H$ is $b$, which is  smooth
 in the strictly convex case and linear-on-faces in the polygonal case.
 In fact, a little reflection shows that in Case 1 we can take $K$ to be
 the maximum of the slopes attained by taking $X_{1}, X_{2}, X_{3}$ to be
   vertices of the polygon.
     So we can apply the result of \cite{kn:GT} quoted above, and deduce that the
     derivative of $H$ is bounded by the constant  throughout $\Omega$.
     Thus the vector field $w$ satisfies an $L^{\infty}$ bound over $\Omega$
     and hence the same is true of $v$.
         To sum up we have
\begin{thm}
Suppose $u\in \cS_{\Omega,\sigma}$ is a solution of Abreu's equation  where $A$ is
a constant and the dimension $n$ is $2$. Then there is a constant $K$,
depending tamely on $\Omega, \sigma$,  such that
$\vert v\vert \leq K$ throughout $\Omega$.
\end{thm}

\subsection{Integrating bounds on $v$}

In this subsection we study the implications for $L=\log \det u_{ij}$ 
of an $L^{\infty}$ bound on
the vector field $v^{i}= -u^{ij}_{j}$ associated to a convex function $u$
on a convex domain $D\subset \bR^{n}$. We prove two results, under different
hypotheses. For the first, recall that the derivative $\nabla u$ is a diffeomorphism
from $D$ to its image $(\nabla u)(D)\subset \bR^{n}$.
\begin{lem}
  Suppose that the image $(\nabla u)(D)$ is convex and $\vert v \vert\leq
  K$ on $D$. Then for any $x,y\in D$
  $$  \vert L(x) - L(y) \vert \leq K \vert (\nabla u)(x)- (\nabla u)(y)
  \vert . $$
  \end{lem}
   To see this, let $\xi_{i}$ be the standard
    Euclidean co-ordinates on the image $(\nabla
   u)(D)$. (More invariantly, $\nabla u $ should be thought of as mapping
    $D$ to the dual space 
   $\left(\bR^{n}
   \right)^{*}$ and $\xi_{i}$ are just the dual co-ordinates of $x^{i}$.)
   In traditional notation, $\xi_{i}=u_{i}$ so
   $$  \frac{\partial \xi_{i}}{\partial x^{j}} = u_{ij}\ , \frac{\partial
   x^{j}}{\partial \xi_{i}} = u^{ij}. $$
   Now 
      $  v^{i} = u^{ij}L_{j}$ (see (8)) but
      by the chain rule 
      $$  u^{ij}L_{j} = \frac{\partial L}{\partial \xi_{i}}. $$
      In other words, if we define a function $L^{*}$ on $(\nabla u)(D)$
      to be the composite $L\circ (\nabla u)^{-1}$ then $v$ is the derivative
      of $L^{*}$ and hence $\vert \nabla L^{*} \vert \leq K$. If $(\nabla
      u)(D)$ is convex then $\nabla u(x), \nabla u(y)$  can be joined by
      a line segment in $\nabla u(D)$ and the result follows immediately
      by integrating the derivative bound along the segment.
      
    The next results applies to solutions of Abreu's equation, but avoids
    the convexity hypothesis on the image under $\nabla u$. We write
    $ \Av_{D}(L)$ for the mean value
    $$ \Av_{D}(L)= (\Vol (D))^{-1} \int_{D} L. $$  
\begin{thm}
Suppose that $D\subset \bR^{n}$
 is a bounded convex domain with smooth boundary and
  $x_{0}\in D$ is a base point.
 Let $R$ be the maximal distance
 from $x_{0}$ to a point of $\partial D$.
 Let $u$ be a solution of Abreu's equation $S(u)=A$ in $D$, smooth up to
 the boundary and normalised at $x_{0}$. 
Then
  $$ \vert L(x_{0}) - \Av_{D}(L) \vert \leq C \int_{\partial D} u d\nu,
  $$
  where $d\nu$ is the standard Riemannian volume form on $\partial D$ and
  $$ C= (n \Vol(D))^{-1}\left( n \sup_{D} \vert v \vert + R  \sup_{D} 
  \vert A \vert \right).
   $$
  \end{thm}
  To prove this, we can suppose that $x_{0}$ is the origin and let
  $(r,\utheta)$ denote standard \lq\lq polar-coordinates'' on 
  $\bR^{n}\setminus\{0\}$
  (so $\utheta\in S^{n-1}$). Suppose that $\partial D$ is given in these
  co-ordinates by an equation $r=\rho(\utheta)$.
  Let $f$ be the function 
  $$ f= \left( 1- \frac{\rho(\utheta)^{n}}{r^{n}}\right)$$
  on $\bR^{n}\setminus\{0\}$ and define a vector field $\zeta$ on 
  $\bR^{n}\setminus\{0\}$  by
$$  \zeta^{i} = f  x^{i}. $$
Then, away from the origin, $$ \zeta^{i}_{i}= 
x^{i}_{i}(1-\frac{\rho^{n}}{r^{n}}) - x^{i} \frac{\partial}{\partial
x^{i}} \frac{\rho^{n}}{r^{n}}= n(1-\frac{\rho^{n}}{r^{n}}) - r 
\frac{\partial}{\partial
r} \frac{\rho^{n}}{r^{n}}= n, $$
and $\zeta$ vanishes on the boundary of $D$. It is clear that
 $\zeta$  satisfies
the distributional equation  
$$ \zeta^{i}_{i} = n- n \Vol(D)\  \delta_{0}, $$
where $\delta_{0}$ is the delta-function at the origin. Applying this distribution
 to
the function $L$ we get
$$  n (\Vol(D) L(0) - \int_{D} L) = \int_{D} \zeta^{i}L_{i}. $$

Now recall that $h=u- u_{i} x^{i}$ satisfies $h_{i} = -u_{ij}x^{j}$. Then
$L_{i}= u_{ij} v^{j}$, and
$$    x^{i}L_{i}= u_{ij} x^{i} v^{j}= - h_{j} v^{j}= 
- \left( hv^{j}\right)_{j} + h v^{j}_{j}. $$
Using Abreu's equation in the form $v^{j}_{j}=A$ we have
$$     x^{i}L_{i}= -\left( hv^{j}\right)_{j} +Ah. $$
Now  
$$\int_{D} \zeta^{i}L_{i} = \int_{D} f x^{i}L_{i}=
 -\int_{D} f (hv^{j})_{j}
 + \int_{D} Afh .$$
 We can integrate by parts on the first term to get
 $$\int_{D} \zeta^{i} L_{i} = \int_{D} (h f_{j} v^{j} + Afh). $$
 There is no boundary term from $\partial \Omega$ 
 since $f$ vanishes on the boundary. There is also no extra term caused by the
 singularity of $f$ at the origin since $h$ vanishes to second order there.
 To sum up, we have the identity, 
 \begin{equation}  n \Vol(D)\left( L(0)-  \Av_{D}(L) \right) = (I) + (II), \label{eq:intid}\end{equation}
 where $$(I)= \int_{D} h f_{j} v^{j} \ \ \ \ (II)= \int_{D}A f h. $$

 We next estimate the size of the integrals $(I)$ and $(II)$. We have
 $$  \frac{ \partial f}{\partial r} = -n\frac{\rho^{n}}{r^{n+1}}\ \ , \
  \frac{1}{r}
 \frac{\partial f}{\partial \utheta} = n \frac{\rho^{n-1} }{r^{n+1}} \frac{\partial
 \rho}{\partial \utheta}, $$
 so $$ \vert \nabla f\vert = \frac{n \rho^{n-1}}{r^{n+1}} W(\utheta) $$
 
  where $$W(\theta) = \sqrt{ \rho^{2} + \vert \rho_{\utheta}^{2}}. $$
   Now let $K=\sup_{D} \vert v \vert$ so
   we have
  $$  \vert f_{j}v^{j} \vert \leq \frac{nK \rho^{n-1}}{r^{n+1}} W(\utheta). $$
  Thus 
  $$ \vert(I)\vert\leq nK \int_{D} W(\utheta) \rho(\utheta)^{n-1}
   \frac{h}{r^{n+1}}= 
  2K\int_{ D} \rho^{n-1}
  W(\utheta)
  \frac{h}{ r^{2} }  dr d\utheta, $$
  in an obvious notation. 
  Now consider this integral along a fixed ray $\utheta={\rm constant}$. In polar
  co-ordinates we can write
  $$ h= u - r \frac{\partial u}{\partial r}, $$
  so $$ \frac{\partial}{\partial r} (r^{-1} u)= -r^{-2} u + r^{-1}
  \frac{\partial u}{\partial r} = - r^{-2} h. $$
  So $$ \int_{\epsilon}^{\rho(\utheta)} \frac{h(r)}{r^{2}} dr= \rho^{-1}
   u(\rho(\utheta),
  \utheta)- \epsilon^{-1} u(\epsilon,\utheta), $$
  and the term from the lower limit tends to zero with $\epsilon$ since
  the derivative of $u$ vanishes at the origin. So we have
  \begin{equation} \vert (I)\vert\leq nK \int_{\partial D} \rho(\utheta)^{n-2} W(\theta) u
  d\utheta. \label{eq:Ibound} \end{equation}
  
  For the other integral we have
        $ - r u_{r} \leq h \leq 0,$ and $f\leq 0$ so
        \begin{equation} 0\leq (II)\leq \sup \vert A \vert
        \int_{D} \frac{\rho^{n}}{r^{n}} ru_{r} r^{n-1} dr d\utheta=
         \sup \vert A \vert \int_{\partial D}
          \rho(\utheta)^{n} u d\utheta. \label{eq:IIbound}\end{equation}       
        
        We now transform the integrals over $\partial D$ to the Riemannian
        area form $d\nu$.
        By straightforward calculus
        $$   d\nu = \rho^{n-1} \sqrt{ 1 + 
        \left(\frac{\rho_{\utheta}}{\rho}\right)^{2}} d\utheta = \rho^{n-2}
        W(\utheta) d\utheta.  $$
        So inequality~\ref{eq:Ibound} is just
        $$ \vert (I) \vert \leq nK \int_{\partial D} u d\nu. $$
        Also, $ d\nu \geq \rho^{n-1} d\theta$ so inequality ~\ref{eq:IIbound}
        gives
             $$ \vert (II)\vert \leq \sup \vert A \vert \int_{\partial D}
             u \rho d\nu \leq R \sup \vert A \vert \int_{\partial D} u
             d\nu, $$ 
       and putting these together gives the result stated.

        We now apply these results in the case when $n=2$ and $A$ is a
        constant, using Theorem 2. We obtain
        \begin{thm}
        Let $\Omega\subset \bR^{2}$ and $u\in \cS_{\Omega,\sigma}$ be a normalised
         solution
        of Abreu's equation with $A$ constant. Suppose Condition 1 holds.
         Then
        $$  L  \leq C_{0}+ C_{1} \vert \nabla u \vert \leq  C_{2} d^{-2}, $$
        where $d$ is the distance to the boundary of $\Omega$ and $C_{0},C_{1},
        C_{2}$ 
        depends tamely on $\Omega, \sigma, \lambda$.
        \end{thm}
        We can prove this using the simpler result Lemma 3. Any function $u$
        in $\cS_{\Omega,\sigma}$ is continuous on $\overline{\Omega}$ and so attains
        a minimum value on $\overline{\Omega}$, but it is clear from the
        definition of $\cS_{\Omega,\sigma}$ that this cannot be attained on the
        boundary. Changing $u$ by the addition of any linear function we
        see that for $u\in \cS_{\Omega,\sigma}$ the image $(\nabla u)(\Omega)$
        is the whole of $\bR^{n}$. Let $B$ be a small disc about the base
        point. We know that $\nabla u $ is bounded on $B$ by Corollary
        3, and this
        gives a bound on the area of $\nabla u(B)$ , which is
        $$  \int_{B} \det u_{ij} . $$
        So there is a $c$, depending tamely on the data,
         such that we can find some point in $B$ where $L\leq c$.
          Now the result follows by combining
        Corollary 3, Theorem 2 and Lemma 3.
        
        We can give a slightly different proof of the above result using
        Theorem 3 in place of Lemma 3. In the case when $\Omega\subset \bR^{2}$ is a polygon and $A$ is
        a constant we can also obtain an interior lower bound on $L$ by applying
        Corollary 4, which gives a bound on
        $$ \int_{\Omega} \vert L \vert. $$
        However we will not discuss this in detail, since we will get a
        better lower bound in the next section.

  \section{Applications of the maximum principle}
  \subsection{Lower bound for the determinant}
  
  In this Section we will derive a lower bound for $L=\log \det(u_{ij})$,
  valid in any dimension and for any bounded function $A$. 
  
  \begin{thm}
  Suppose  that
  $u\in \cS_{\Omega,\sigma} $ is a solution of Abreu's equation for some
   bounded function $A$. For any $\alpha\in (0,1)$
   there is a constant $C_{\alpha}$ (depending only on $n$ and $\alpha$)
   such that 
      $$ \det u_{ij} \geq C_{\alpha} d^{-\alpha} (\sup A)^{n} \ 
      {\rm Diam}(\Omega)^{2n+\alpha} $$ throughout $\Omega$.
\end{thm}
Here we recall that $d$ is the distance to the boundary of $\Omega$.
${\rm Diam}\ (\Omega)$ is the diameter of $\Omega$. 

We can compare this result with the {\it asymptotic behaviour} established
in Section 2 (Proposition 2 and (20)). In Case 1
\begin{equation} \det u_{ij} \sim C_{1} d^{-1}, \label{eq:case1} \end{equation}
as $x$ tends to a generic boundary point (i.e on an $(n-1)$-dimensional face
of the boundary) and in Case 2,
\begin{equation}
 \det u_{ij}\sim C_{2} d^{-1} \vert \log d\vert. \label{eq:case2}\end{equation}

To prove the theorem we suppose that we have a function $\psi$ on $\Omega$
 which satisfies
the following conditions.
\begin{enumerate}
\item For each point of $\Omega$ the matrix
$$  M_{ij}= \psi_{ij} - \psi_{i}\psi_{j} $$
is positive definite.
\item $L-\psi$ tends to $+\infty$ on $\partial \Omega$. \end{enumerate}
Certainly such functions exist. For example, we can take
$\psi_{\epsilon}(x) = \epsilon \vert x \vert^{2}$ for small $\epsilon$: the second
condition holds since $L\rightarrow \infty$ on $\partial \Omega$. Given
such a function $\psi$, we write $D=\det (M_{ij})$.

By the second condition there is a point $p$ in $\Omega$ where $L-\psi$
attains its minimum value. At this point we have
$$   L_{i}=\psi_{i} \ \ \ (L_{ij}- \psi_{ij})\geq 0. $$
We take Abreu's equation in the form (Equation~\ref{eq:Ab4})
$$   u^{ij} ( L_{ij}- L_{i}L_{j}) = A. $$
At the point $p$ we have
$$   u^{ij}L_{ij}= A+ u^{ij} L_{i}L_{j} = A+ u^{ij} \psi_{i}\psi_{j}, $$
and $$u^{ij} L_{ij} \geq u^{ij} \psi_{ij}.$$
Thus, at the point $p$,
$$   u^{ij} (\psi_{ij}-\psi_{i}\psi_{j}) =u^{ij}M_{ij}\leq A. $$
So $ u^{ij}M_{ij}\leq \overline{A}$, where
$\overline{A}=\sup_{\Omega} A$. We now use the standard inequality for
positive definite matrices,
\begin{equation}   \left( \det(u^{ij}) \det(M_{ij}) \right)^{1/n}\leq
 n^{-1} u^{ij}M_{ij}. \label{eq:AGM}\end{equation}
 (This is just the arithmetic-geometric mean inequality for the relative
 eigenvalues of $(u_{ij}), (M_{ij})$.) Now use the fact that $\det(u^{ij})
 = e^{-L}$ to get, at the point $p$,
 $$    e^{-L} D \leq \left( \overline{A}/n\right)^{n}, $$
 or $$ L(p)\geq  \log D(p) -c, $$
 where  $c= n\log(\overline{A}/n)$. Then for any other point, $q$, in
 $\Omega$ we have $  (L-\psi)(q) \geq (L-\psi)(p)$ so
 $$  L(q) \geq L(p) + \psi(q) - \psi(p) \geq \log D(p) +\psi(q)-\psi(p)-c,
 $$
 or in other words
 \begin{equation} L(q) \geq \psi(q) + C_{\psi}, \label{eq:lowbound}\end{equation}
 where $$ C_{\psi} = \inf_{\Omega} (\log D - \psi)- c. $$
 (Of course, at this stage, $C_{\psi}$ could be $-\infty$, in which case
 Equation~\ref{eq:lowbound} is vacuous.)

 By taking the function $\psi_{\epsilon}$, say, we immediately get a lower bound
  $L\geq\ {\rm const.}$ over $\Omega$. (This is all we need for our main
  application below.) To prove Theorem 5 we need to make a
 more careful choice of the comparison function $\psi$. In fact the optimal
 choice
 of this comparison function leads to an interesting Monge-Amp\`ere differential
 inequality which we will digress to explain. We consider, for a fixed point $q$ the set of
 functions $\psi$ satisfying conditions (1) and (2) and with $C_{\psi}>-\infty$.
 The optimal bound we can get from the argument above is given by the supremum
 over this set of functions $\psi$ of $\psi(q)+ C_{\psi}$. Changing $\psi$
 by the addition of a constant, we may suppose that $C_{\psi}=0$. In other
 words we have
 $  L(q) \geq \lambda(q)- c$, where
 $$ \lambda(q) = \sup_{\psi\in X} \psi(q) $$
 and $X$ is the set of functions which satisfy (1),(2) and in addition
 \begin{equation}  \log D \leq \psi . \label{eq:bound}\end{equation}
 This becomes more familiar if we write $R=-\exp(-\psi)$. Then
 $$   R_{ij} = e^{-\psi}( \psi_{ij}- \psi_{i}\psi_{j})= e^{-\psi}M_{ij}.
 $$  So condition (1) is simply requiring that $R$ be a convex function.
 We have $D= e^{n\psi} \det(R_{ij})$ so Equation~\ref{eq:bound} becomes
 $$     \det(R_{ij}) \geq (-R)^{n-1}.  $$
 
 We can sum up the discussion in the following way. Given a domain
 $\Omega$,
  let ${\cal R}_{\Omega}$ denote the set of negative 
 convex functions $R$ on
 $\Omega$ with 
 $$\lim_{x\rightarrow \partial \Omega} \frac{R(x)}{d(x)}= -\infty$$
and with $$ \det (R_{ij}) \geq (-R)^{n-1}. $$
Define a function $\rho_{\Omega}$ by
$$  \rho_{\Omega}(x) = -\sup_{R\in {\cal R}_{\Omega}} R(x). $$
Then we have
 \begin{prop}
 In either Case 1 or Case 2, a solution $u\in \cS_{\Omega,\sigma}$ 
 of Abreu's equation with boundary data 
 satisfies
  $$ \det u_{ij} \geq  \rho_{\Omega}^{-1} \left( \frac{\sup_{\Omega} A}{n}\right)^{n}.
  $$ in $\Omega$.
  \end{prop}
  
  This follows from the argument above and the asymptotic behaviour 
  ~\ref{eq:case1},
~\ref{eq:case2}  
  (In Case 2 we could strengthen the statement a little, taking account
  of the $\log d$ term in ~\ref{eq:case2}.)
 
  We return from this digression to complete the proof of Theorem 5.
  For this we take co-ordinates $(x_{1}, \dots ,x_{n-1},y)$ on $\bR^{n}$
  and consider the function
  $$ r(\ux, y)=  y^{\alpha} ( \frac{b}{2} \ux^{2} - 1 ), $$
  on a cylinder $Z=\{(\ux,y): \vert \ux\vert\leq 1, 0<y<1 \}$.
   Here $\alpha \in (0,1)$ is given, as in the statement of the Theorem,
    and $b>0$ will be specified later. We will choose $b$ with $b< 2$
  so that $r$ is negative on $Z$.  
   Straightforward calculation shows that $r$ is convex on $Z$ provided
   $(1-\alpha) (1-\frac{b}{2}) \geq  \alpha b$. This is equivalent to
   $$  \frac{2-b}{2+b} > \alpha, $$
   so, whatever $\alpha$ is given, we can choose $b$ so small that the
   condition holds. Then 
   $$ \frac{\det (r_{ij})}{(-r)^{n-1}} = \alpha y^{\alpha-2} \left( (1-\alpha) - \alpha 
   \frac{\ux^{2}}{1-\frac{b}{2} \ux^{2}}
   \right) (\frac{b}{1-\frac{b}{2}\ux^{2}})^{n-1}, $$
   which is bounded below on $Z$ so  
   $\det (r_{ij}) \geq C (-r)^{n-1}$ throughout $Z$  for some fixed positive
   $C$. Writing $R= C^{-1} r$ we get a function with
   $\det R_{ij} \geq (-R)^{n-1}$ on $Z$ and with $-R=O(y^{\alpha})$ as
   $y\rightarrow 0$.  
   
   Now one checks first that the inequality to be proved is stable under
   rescaling of the domain $\Omega$. Then given a point $p$ in our domain
    $\Omega$ we can obviously suppose without loss of generality 
   that $\Omega$ lies in the set $Z$ (for some suitable $K$ depending on
   $\Omega$), that the origin lies in $\partial \Omega$ 
    and that the minimum distance to the boundary of $\Omega$
   from $p$ is achieved at the origin. Then the function $R$ constructed
   above (or, more precisely, its restriction to $\Omega$) lies in the
   set ${\cal R}_{\Omega}$, since $R/d$ tends to $-\infty$ on $\partial
   \Omega$,   and the result follows from Proposition
   5.

   It should be possible to sharpen this bound in various ways, and 
    this is  related
   to  the Monge-Amp\`ere Dirichlet problem  for convex, negative, functions $R$ on
   $\Omega$:
   $$   \det R_{ij} = (-R)^{n-1} \ \  R\vert_{\partial\Omega}=0.
   $$

\subsection{An upper bound on the determinant}

In this subsection we will modify the method of Trudinger and Wang in \cite{kn:TW1},
\cite{kn:TW2} to obtain
an upper bound on the determinant of $(u_{ij})$. The argument applies in
any dimension but the result we obtain requires additional information
on \lq\lq modulus of convexity'' of $u$ for its application, see (5.1)
below.

Consider a bounded domain with smooth boundary $D$ in $\bR^{n}$ and a
 smooth convex function $u$ on
$D$ satisfying $S(u)=A$, with $u<0$ in $D$ and $u\rightarrow 0$ on
 $\partial D$. 
We suppose that $u$ is smooth up to the boundary of $D$. We consider
Euclidean metrics $g_{ij}$ on $\bR^{n}$ with determinant $1$ and for
each such metric let 
$$   C_{g}= \max_{D} g^{ij} u_{i} u_{j}. $$
Now let 
$$ C= \min_{g} C_{g}. $$
This defines an invariant  $C$ of the function $u$ on $D$. Another way
of expressing the definition,  is that $\omega_{n}C^{n}$ is least volume
of an ellipsoid containing the image of $D$ under the map $\nabla u:D\rightarrow
\bR^{n}$, where $\omega_{n}$ is the volume of the unit ball in $\bR^{n}$.

\begin{thm} In this situation 
$$  \left( \det(u_{ij})\right)^{1/n}\leq
 \left(\frac{5}{2} +\frac{aM}{2n}\right) e C \ (-u)^{-1}, $$
  
in $D$, where
$$a= \max (0, - \min_{D} A)\ \ , \ \ M= \max_{D}(-u). $$

\end{thm}

To prove this we let $g$ be the metric with $C_{g}=C$ and consider the
function
$$   f=-L - n \log(-u) - \alpha g^{ij}u_{i}u_{j} $$
on $U$, where $\alpha$ is a positive constant to be fixed later.
 Thus $f$ tends to $+\infty$ on $\partial D$ and there is a point
$p$ in $D$ where $f$ attains its minimum. At the minumum 
$f_{i}=0$ and $u^{ij} f_{ij} \geq 0$. The first of these gives
\begin{equation}  L_{i} = -\frac{n u_{i}}{u} - 2\alpha g^{pg} u_{p} u_{qi}. \label{eq:deriv}
\end{equation}
The second gives
\begin{equation} 0\leq - u^{ij} L_{ij} - \frac{n^{2}}{u} + \frac{n}{u^{2}} u^{ij}u_{i}u_{j}
- 2\alpha ( u^{ij}g^{pq} u_{pi} u_{qj} + u^{ij}g^{pq} u_{p} u_{qij}),\label{eq:2deriv}\end{equation}
where we have used the fact that $u^{ij}u_{ij}= n$.
The crucial step now is to observe that
$$u^{ij} g^{pq} u_{pi} u_{qj} = g^{pq}u_{pq}, $$
the ordinary Euclidean Laplacian of $u$, in the metric $g$.
We next use the form Equation~\ref{eq:Ab4} of Abreu's equation 
to see that 
\begin{equation}  u^{ij}L_{ij}\geq u^{ij} L_{i}L_{j}+ \uA, 
\label{eq:Abineq}\end{equation}
where $\uA= \min_{D} S(u)$.
Now Equation~\ref{eq:deriv} gives, at the point $p$,
$$ u^{ij}L_{i} L_{j} = u^{ij} \left(  \frac{n^{2}}{u^{2}} u_{i} u_{j}
+ 4\alpha^{2} g^{pg}g^{rs} u_{p} u_{r} u_{qi} u_{sj} 
+ 4\frac{\alpha n}{u} g^{pq}u_{i} u_{p} u_{qj}\right), $$
which can be written in the simpler form
$$ u^{ij} L_{i} L_{j}= \frac{n^{2}}{u^{2}} u^{ij}u_{i}u_{j} + 4\alpha^{2}
g^{pq}  g^{rs} u_{pr}u_{q} u_{s} + 4\frac{\alpha n}{u} g^{pq} u_{p} u_{q}.
$$
Our inequalities~\ref{eq:2deriv},~\ref{eq:Abineq} give
\begin{eqnarray}   0\leq &-& \uA-\frac{n^{2}}{u^{2}} (u^{ij} u_{i} u_{j}) - 4\alpha^{2}(g^{pq}
u_{rs} u_{p} u_{q}) - \frac{4\alpha n}{u} g^{pq}u_{p}u_{q} - \frac{n^{2}}{u}
+\nonumber\\&+& \frac{n}{u^{2}} (u^{ij} u_{i} u_{j}) - 2\alpha g^{pq}u_{pg} - 2\alpha
u^{ij}g^{pq} u_{pij} u_{q}. \label{eq:massive} \end{eqnarray}

Now recall that
$L_{q} = u^{jk} u_{qjk}$, so Equation~\ref{eq:deriv} gives
$$ g^{pq} u^{ij} u_{pij} u_{q} =  g^{pq}  L_{p} u_{q} = 
-  \frac{n}{u} g^{pq} u_{p} u_{q} -
 2\alpha u_{pq} u_{r} u_{s} g^{pr} g^{qs}$$ 

This means that three of the terms in (40) cancel. There are also two terms
involving the expression $u^{ij} u_{i} u_{j}$ which we can combine 
 to get
\begin{equation}  0\leq -\uA - \frac{n^{2}}{u} + \frac{(n-n^{2})}{u^{2}} (u^{ij}u_{i}u_{j})
-2\frac{\alpha n}{u} g^{pq} u_{p} u_{q} - 2\alpha g^{pq} u_{pq}. \label{eq:simpler}\end{equation}
Now use the fact that $n-n^{2}\leq 0$ and
$$   g^{pq} u_{p} u_{q} \leq C, $$
by the definition of $C$. Re-arranging, we obtain
$$  \frac{g^{pq} u_{pq}}{n} \leq \left( \frac{n}{2\alpha} + C)\right) \left(
\frac{-1}{u}\right) - \frac{\uA}{2\alpha n}, $$
where $\uA=\min A$. By the definition of $a$ and $M$ we have
$$ \frac{g^{pq}u_{pq}}{n} \leq \left(\frac{-1}{u}\right) \left( \frac{n}{2\alpha}
+ C + \frac{aM}{2\alpha} \right). $$
Now, since $g^{pq}$ has determinant $1$ we have (as in Equation~\ref{eq:AGM})
$$ e^{L}= \det(u_{ij}) \leq \left(\frac{ g^{pq}u_{pq}}{n}\right)^{n}, $$
so, at the point $p$,
$$ L\leq -n \log(-u)+ n\log P, $$
where $P=\frac{n+aM}{2\alpha} + C$. 
Going back to the definition of $f$
we obtain
$$ f(p) \geq -\alpha g^{pq} u_{p} u_{q} - n\log P\geq -\alpha C -n\log P.$$
So at any other point  of $D$, we also have$ f\geq - \alpha C- n\log P$
which gives
  $$ \det(u_{ij})^{1/n}\leq  \kappa C (-u)^{-1}, $$
  with $ \kappa = \exp( \frac{\alpha C}{n} + \log P). $
  Optimising the choice of $\alpha$ one finds that the least value of $\kappa$
  is
  $$ \kappa_{{\rm min}} = C \left( \sqrt{m + \frac{m^{2}}{4}} + \frac{m}{2}
  + 1\right) \exp\left( \sqrt{m+\frac{m^{2}}{4}} - \frac{m}{2}\right), $$
  where $m= \frac{1}{2} + \frac{aM}{2n}$. 
  Now $\kappa_{\rm min}\leq Ce (m+2)$ since
   $\sqrt{m + \frac{m^{2}}{4}}\leq \frac{m}{2} + 1$, and this gives the
   inequality stated in the Theorem ({\it i.e.} the statement is not quite the
   optimal result: notice that if $A\geq 0$, so $a=0$, one gets $\kappa_{\rm
   min}= 2 Ce^{1/2}$.)

  \section{Estimates for higher derivatives}
\subsection{The modulus of convexity and the proof of Theorem 1} 

Suppose we have upper and lower bounds on $\det(u_{ij})$ in the interior of
$\Omega$ for a solution $u$ of Abreu's equation, for any smooth function
$A$ and in any dimension. We can then obtain estimates on all derivatives,
following the discussion of Trudinger and Wang \cite{kn:TW1}, \cite{kn:TW2}, in terms of a \lq\lq
modulus of convexity'' of $u$. This uses the results of Caffarelli \cite{kn:Caf}
and
Caffarelli and Guti\'errez \cite{kn:CafGut} and the method exploits  the form 
Equation~\ref{eq:Ab5}, $Q(F)=-A$, of Abreu's equation. 

We introduce some
notation which we will use more extensively in (5.4) below. For  a smooth strictly convex function $u$ on any open set $U\subset
\bR^{n}$ and a point $x\in U$ we let $\lambda_{x}$ be the affine linear
function defining the supporting hyperplane of $u$ at $x$ ({\it i.e.}\
 $u-\lambda_{x}$ vanishes to first order at $x$). Then we define the function
 $H_{x}$ on $U$ by
 $$   H_{x}(y) = u(y) - \lambda_{x} (y). $$
 (So $H_{x}$ is the normalisation of $u$ at $x$, in our previous terminology.)
Thus $H_{x}(y)\geq 0$ with equality if and only if $x=y$. We can think
of the function $H_{x}(y)$ of two variables $x,y$ as a kind of \lq\lq distance
function'' on $U$, although it need not satisfy the axioms of a metric.
For a subset $S\subset U$ we put $H_{x}(S) = \inf_{y\in S} H_{x}(y)$.

Now return to the case where $u$ is a convex function on $\Omega$ and $K\subset \subset
K^{+}$
 are compact convex subsets of $\Omega$ with  $0<\lambda\leq \det u_{ij}\leq \Lambda$
on $K^{+}$. Let
 $$ H(K,K^{+}) = \min_{x\in K} H_{x}(\partial K^{+}). $$  Caffarelli and 
 Guti\'errez
 prove that there is an $\alpha\in (0,1)$ such that
 for any solution $f$ of the linear equation $Q(f)=-A$ and $x,y $ in $K$
 $$ \vert f(y)- f(x) \vert \leq C H(K,K^{+})^{-\alpha} \vert x-y \vert^{\alpha}$$
 where $C$ depends only on $A$, the supremum of $\vert \nabla u\vert$ over
 $K^{+}$ and the upper and lower bounds $\Lambda,\lambda$ of
  $\det(u_{ij})$ over $K$. 
[The results of \cite{kn:CafGut} are stated for the case when $A=0$ but, according to
Trudinger and Wang (\cite{kn:TW2}, discussion following Lemma 2.4),
 the arguments go over to the inhomogeneous equation.]
Thus in our situation, taking $f=F=\det(u_{ij})^{-1}$,
 if we have a positive lower bound on $H(K,K^{+})$ we get a Holder
estimate on $F$ and hence on $\det(u_{ij})$. Then the results of Caffarelli
in \cite{kn:Caf}
give  $C^{2,\alpha}$ bounds on $u$ over $K$, depending again on the $H(K,K^{+})$.
This gives $C^{\alpha}$ control of the co-efficients of the linearised
operator $Q$ and we can apply the Schauder estimates (\cite{kn:GT} Theorem
6.2) to get $C^{2,\alpha}$
control of $F$, and so on. 

Further, if we have a lower bound on the \lq\lq modulus of convexity''
$H(K,K^{+})$ we can apply Theorem 6, using the interior bound Corollary
3 on the first
derivatives of $u$, to obtain an upper bound on $\det(u_{ij})$ over $K$.
 Given a point $x$ in $K$ we set 
$$\tilde{u} = H_{x}- \frac{1}{2} H(K,K^{+})$$ and
$$ D= \{ y: \tilde{u}(y) \leq 0\}, $$
so $D\subset K^{+}$ and we can apply Theorem  6 to $\tilde{u}$ to get an
upper bound on $\det(u_{ij})$ at $x$. Combining this with the lower bound
from (4.1) (and replacing $K^{+}$ by $K$) we can then feed into the
 the preceding argument. Thus
the obstacle to proving a result like Theorem 1 in general dimensions $n$ is the
need to control the $H(K,K^{+})$, for interior subsets $K\subset \subset
K^{+} \subset \Omega$. In dimension $2$ we can again argue in parallel
with Trudinger and Wang in \cite{kn:TW2}. The result of (4.1) gives a lower bound $\det(u_{ij})\geq
\lambda_{K^{+}}>0$, over $K^{+}$, 
for a solution satisfying our boundary conditions. Then
 a result of Heinz \cite{kn:Heinz}, implies that in two dimensions,
    $$ H(K,K^{+})\geq C \lambda_{K^{+}}^{1/2}, $$
    where the constant $C$ depends on $K^{+}$, the distance from $K$ to the boundary
    of $K^{+}$ and $\sup_{K^{+}} \vert \nabla u \vert$.
     Putting all these facts from the literature together, we arrive
    at a proof of our main Theorem 1.
    
      As we explained in the Introduction,
     we will give below an alternative proof, in the case when $A$ is constant
    and $\Omega$ is a polytope, which avoids the sophisticated analysis
    of the Caffarelli-Guti\'errez theory but employs some of their  basic
    tools. The modulus of convexity will also enter in a crucial
    way in our argument. We close this subsection with two further remarks
    \begin{itemize}
     \item In the case when $A$ is a constant and $\Omega$ is a polygon
     in $\bR^{2}$ we can combine the lower bound Theorem 5 and the upper bound
     Theorem 4, to obtain
     $$ \vert \nabla u \vert \geq C d^{-\alpha} $$
     near the boundary of $\Omega$, for some $C,\alpha>0$. This  gives easily  a lower bound on
     $H(K,K^{+})$ for  suitable $K, K^{+}$ and means we can avoid appealing
     to the result of Heinz in this case. 
     \item In two dimensions the main result of Heinz in \cite{kn:Heinz}
      gives a $C^{1,\beta}$
     bound on $u$, in terms of upper and lower bounds for $\det(u_{ij})$.
     In the case when $A$ is a constant we can use Lemma 3, and the bound
     on the vector field $v$ (Theorem 2) to obtain a $C^{,\beta}$ bound
     on  $\det(u_{ij})$.
     This gives an alternative path, avoiding appeal to \cite{kn:CafGut}, but feeding
     into \cite{kn:Caf},  in this case.   
\end{itemize}

\

The strategy for our alternative proof is as follows. In (5.2) we introduce
two tensors $F,G$  depending on the 4th. order derivatives of the function
$u$ and related to the Riemann curvature and Ricci tensors of a certain Riemannian
metric. We show that our boundary conditions fix the natural $L^{2}$ norms
of these tensors. Then we show in (5.2) that under suitable conditions
the $L^{\infty}$ norm of $G$ (or $F$) controls the second derivatives of
the function $u$, while in (5.3) we show that under suitable hypotheses,
including control of the second derivatives of $u$, the $L^{\infty}$ norm
of $G$ is  controlled by the $L^{2}$ norm of $F$. Putting together
these three ingredients we complete the proof in (5.4), using a scaling
argument and some of the basic geometrical results of Caffarelli and Guti\'errez
on the sections of a convex function.    

\subsection{Curvature identities and $L^{2}$ bounds}

For a convex function $u$ on an open set $U$ in $\bR^{n}$ we define a 4-index
 tensor
by
   $$   F^{ab}_{kl} = -u^{ab}_{kl}. $$
 We can raise and lower indices in the usual way, using the metric $u_{ij}$,
 setting
 $$   F^{abcd} = u^{ck} u^{dl} F^{ab}_{kl}\ , F_{ijkl} = u_{ia} u_{jb}
 F^{ab}_{kl}. $$
 \begin{lem}
 The tensor $F_{ijkl} $ is symmetric in pairs of indices
 $$  F_{ijkl}= F_{klij}. $$
 \end{lem}
 In fact, calculation gives,
 $$  F_{ijkl} =u_{ijkl}- u^{\lambda q}( u_{kjq }u_{il\lambda} + 
 u_{ik\lambda}u_{jlq})$$
 which makes the symmetry apparent.
 
We now introduce a Riemannian metric, due to Guillemin \cite{kn:Guil}, on the $2n$-dimensional
manifold $U\times \bR^{n}$, with co-ordinates $\eta_{i}$ in the second
factor;
\begin{equation}   g= u_{ij} dx^{i} dx^{j} + u^{ij} d\eta_{i} d\eta_{j}. \label{eq:metric}
\end{equation}
This is in fact a Kahler metric: complex co-ordinates and a Kahler potential
 are furnished
by the Legendre transfrom construction. If $\xi_{i} = u_{i}$ are the
usual transformed co-ordinates we set
$z_{i} = \xi_{i} + \sqrt{-1}\ \eta_{i}$
so
\begin{equation}  dz_{i} = u_{ij} dx^{j} + \sqrt{-1}\ d\eta_{i}. \label{eq:oneforms}\end{equation}
\begin{lem}
The curvature tensor of $g$ is
$$   -F^{ijkl} dz_{i} d\oz_{k} \otimes dz_{j} d\oz_{l}. $$
\end{lem}

Of course we can use (43) to express the curvature tensor entirely in terms
of products of the $dx^{i}$ and $d\eta_{j}$, avoiding the Legendre transform.

The proof of the Lemma is mainly a matter of notation. We use the standard fact
that if the Kahler metric is expressed by a (Hermitian) matrix-valued function
$H$, in complex co-ordinates, then the curvature, viewed as a matrix of
$2$-forms, is $-\db(H (\partial(  H^{-1})))$.  In our case the matrix $H$ has
entries
$$  H_{\lambda \mu} = 
\langle \frac{\partial}{\partial z_{\lambda}}, \frac{\partial}{\partial
z_{\mu}} \rangle = u^{\lambda \mu}. $$
so $\Gamma= -H(\partial H^{-1})$ is the matrix of $1$-forms with $(\lambda \mu)$
entry 
$$  -u^{\lambda \mu} \frac{\partial u_{\mu \nu}}{\partial \xi_{\alpha}}
 dz_{\alpha}. $$
Now use the fact that $\frac{\partial}{\partial \xi_{\alpha}}= u^{\alpha
i} \frac{\partial}{\partial x^{i}}$ to write this as
$$ \Gamma^{\lambda}_{\mu}=- u^{\lambda \nu} u^{\alpha i}
 u_{\mu \nu i} dz^{\alpha},
$$
and this is just
$$ \Gamma^{\lambda}_{\mu}= u^{\lambda \alpha}_{\mu} dz_{\alpha}. $$

So the curvature, as a matrix of $2$-forms, has $\lambda,\mu$ entry
$$  
 \frac{\partial}{\partial \xi_{\beta}}\left( u^{\lambda \alpha}_{\mu}\right)
 d\oz_{\beta}dz_{\alpha}. $$
 Expressing the derivatives in terms of the $x^{i}$ variables again, this
 just means that the $(1,3)$ curvature tensor is
 $$     u^{\beta j} u^{\lambda \alpha}_{\mu j} dz_{\lambda}\otimes \frac{\partial}{\partial
 z_{\mu}} \otimes d\oz_{\beta} dz_{\alpha}, $$
 and lowering an index gives the stated formula for the $(0,4)$ curvature
 tensor.

 Of course we can also obtain Lemma 3 via Lemma 4 and the usual symmetries
 of the curvature tensor, but we preferred to give the direct calculation.
 In fact in what follows we will make little explicit
  use of the Riemannian metric $g$, and derive most of the formulae we need
  directly.

  Now define a 2-tensor $G$ by contracting $F$:
  \begin{equation}   G^{i}_{k} = F^{ij}_{kj}, \label{eq:Ricci}\end{equation}
  and likewise $G^{ik}, G_{ik}$. If follows from Lemma 3 that the latter are
  symmetric tensors, and contracting in Lemma 4 shows that $G$ is essentially equivalent
  to the Ricci tensor of the metric $g$. In terms of our vector field $v$,
  $$  G^{i}_{k} = v^{i}_{k}. $$
  A further contraction yields the scalar invariant $S= G^{i}_{i}$ of $u$
  which is of course the term appearing in Abreu's equation and which corresponds
  to the scalar curvature of $g$.
  
  It is well-known in Kahler geometry that on a compact Kahler manifold
  the $L^{2}$-norms of the scalar curvature, the Ricci curvature and the
  full curvature tensor give essentially equivalent functionals on the
  metrics in a given Kahler class: they are related by topological invariants
  of the data \cite{kn:Cal}. In our setting we do not necessarily have a compact Kahler
  manifold available so we will develop the corresponding theory directly.
  
  On account of the symmetry of Lemma 3, the standard square-norm of the tensor $F$,
  using the metric $u_{ij}$ is 
  $$  \vert F\vert^{2} = F^{ij}_{kl} F^{kl}_{ij}, $$
  and this is the same (up to a numerical factor) as the standard square-norm
  of the curvature tensor of $g$. Similarly 
  $$  \vert G\vert^{2} = G^{i}_{j} G^{j}_{i}, $$
  is essentially the square-norm of the Ricci tensor. We consider a 1-parameter
  family of functions $u(t)$ with
  $$  \frac{d}{dt} u \vert_{t=0} = \epsilon, $$
  so $\epsilon $ is a function on $U$. We write $E^{ij}$ for the 
  $t$-derivative of $u^{ij}$ at $t=0$, so
  $$      E^{ij} = - u^{ia} \epsilon_{ab} u^{bj}. $$
  
  \begin{prop} The time derivatives at $t=0$ satisfy:
  
   $$\frac{d}{dt} \left( \vert F \vert^{2} - \vert G\vert^{2} \right) = 
   2 Z^{i}_{i}$$
   and $$ \frac{d}{dt} \left( \vert G \vert^{2} - S^{2} \right) = 2 W^{i}_{i}$$
   where 
   $$ Z^{i}= -E^{jl}_{k} F^{ik}_{jl} + E^{ij}_{k} G^{k}_{j}, $$
   and $$ W^{i} = - E^{jk}_{j} G^{i}_{k} + S E^{ji}_{j}. $$ 
   \end{prop}
   The proofs are straightforward calculations.
   
   Now return to the setting of our convex domain $\Omega$ in $\bR^{n}$,
   where we suppose we are in Case 1 with $\Omega$ a polytope.
   Any pair $u_{0}, u_{1}$ of functions in $\cS_{\Omega}$ 
   can be joined by a smooth path $u_{t}$ (for example
   a linear path) in $\cS_{\Omega}$
   For $\delta>0$ let $\Omega_{\delta}$ be an  interior domain 
   with boundary a distance $\delta$ from the boundary
   of $\Omega$. Then the time derivative of
   $$  \frac{1}{2} \int_{\Omega_{\delta}} \vert F \vert^{2} - S^{2} $$
   is the boundary integral
   \begin{equation}   \int_{\partial \Omega_{\delta}} Z^{i} + W^{i}. 
   \label{eq:boundary}\end{equation}
   \begin{lem}
    The boundary integral (45) tends to $0$ with $\delta$, uniformly over
    the parameter $t\in [0,1]$
    \end{lem}
    Given this we obtain 
    \begin{cor} In the case when  $\Omega$ is a polytope there is
    an invariant $\chi(\Omega, \sigma)$ such that for any $u\in \cS_{\Omega,
    \sigma}$
      $$\int_{\Omega} \vert F\vert^{2} - S^{2} = \chi(\Omega, \sigma)
    .$$ \end{cor}
    It is easy to see that $\chi(\Omega, \sigma)$ is a tame function of
    $\Omega, \sigma$. Now if $u\in \cS_{\Omega,\sigma}$ has $S(u)$ constant
    then the value of the constant is fixed by Equation~\ref{eq:Aid} and
    we see that
    $$  \int_{\Omega} \vert F \vert^{2}= \chi(\Omega, \sigma) + \frac{\Vol(\partial\Omega,
    \sigma)^{2}}{\Vol(\Omega)}, $$
    is a tame function of $\Omega, \sigma$.

  To prove Lemma 5 we go back to Proposition 2. Note first  that, in an adapted co-ordinate
  system, the matrix $u^{ij}$ is smooth up to the boundary so
  the same is true for $F^{ij}_{kl}$ and for the time derivative $E^{ij}$.
  (Actually, in Proposition 2 we chose a special adapted co-ordinate sytem to diagonalise
  the Hessian of the function in the variables $x_{i}$ for $i>p$ but it
  is easy to see that the same conclusions hold without this restriction.)
  Thus the vector field $Z,W$ are smooth up to the boundary and we simply
  need to show that the normal components vanish on the boundary. Thus
  we can restrict to a neighbourhood of a point in an $(n-1)$-dimensional
  face of the boundary, i.e. with $p=1$. Now, according to Proposition
  2,
   the entries $u^{1j}$ are all products of $x^{1}$ with smooth functions
   so $u^{1j}_{k}=0$ for all $j$ and $k>1$. Thus 
   \begin{equation}  E^{1j}_{k}= 0\ \ \  k>1 . \label{eq:vanish1} \end{equation}
   The $x^{1}$ derivative of $u^{11}$ is fixed by the boundary measure
   so
    \begin{equation} E^{11}_{1} = 0. \label{eq:vanish2} \end{equation}
   As for the second derivatives we have
   \begin{equation}   F^{1j}_{kl}=-u^{1j}_{kl}= 0\ \ {\rm for}\  k,l>1\label{eq:vanish3}
   \end{equation}
   and \begin{equation} F^{11}_{j1}= -u^{11}_{j1}=0\label{eq:vanish4} \end{equation}
   since the boundary measure is constant. It is now completely straightforward
   to check that (46),(47),(48),(49) imply that all the terms vanish in the sums
   $$  Z^{1}= -E^{jl}_{k} F^{1k}_{jl} + E^{1j}_{k} G^{k}_{j}\ \ ,\ \  W^{1}= -E^{jk}_{j}
   G^{1}_{k}+ S E^{j1}_{j}. $$

  \subsection{Estimates within a section}
  
  We begin by considering a convex function $u$ on a bounded open set,
  with smooth boundary, $D\subset \bR^{n}$. We assume that
  $u$ is  smooth up to the boundary, that $u<0$ in $D$ and that $u$ vanishes
  on $\partial D$. Suppose that, at each point of $D$, we have
  $$ u\geq -c_{0}, \vert G\vert \leq c_{1}, \vert v \vert_{\Euc}\leq c_{2}, \vert \nabla
  u\vert_{\Euc} \leq c_{3}, $$
  for some $c_{0}, c_{1}, c_{2}, c_{3}>0$. Here $G$ denotes the tensor $G_{ij}$
  introduced above and $\vert G\vert$ is the natural norm computed in the
  metric defined by $u$, that is
  $$  \vert G\vert^{2} = G^{i}_{j} G^{j}_{i} = G_{ij} G_{ab}u^{ia}u^{jb}.
  $$On the other hand the quantities $\vert v \vert_{\Euc}, \vert \nabla
  u\vert_{\Euc}$ are the norms of the vector field $v^{i}$ and the derivative
  $u_{i}$ computed with respect to the standard Euclidean metric on $\bR^{n}$.

  \begin{prop}
     There is a constant $K$, depending only on $n, c_{0}, c_{1}, c_{2},
     c_{3}$  such that $(u_{ij})\leq K \vert u \vert^{-1} $ on $D$.      
 \end{prop}
 
 The proof is a straightforward 
 variant of Pogorelov's estimate for solutions of the Monge-Amp\`ere
 equation $\det(u_{ij})=1$, see \cite{kn:Gut} Chapter 4. (It is also similar to the proof
 in (4.2) above.) It obviously suffices to estimate
 the second derivative $u_{11}$ along the co-ordinate axis, and we set
 $$  f= \log (u_{11}) + \frac{u_{1}^{2}}{2}. $$
 Now consider the function $\log (-u) + f$ which tends to $-\infty$ on
 the boundary so has an interior maximum. At this maximum point we have
 \begin{equation}   \frac{u_{i}}{u} + f_{i}= 0 \label{eq:first}\end{equation}
 and $$ P (-\log(u) + f) \leq 0$$
 where $P(\phi )= (u^{ij} \phi_{i})_{j}$. Now
 $$ P \log (u_{11}) = \left( \frac{u^{ij} u_{11i}}{u_{11}}\right)_{j}$$
 which gives
 \begin{equation} P \log (u_{11}) = \frac{u^{ij} u_{11ij}}{u_{11}} + 
 \frac{u^{ij}_{j}u_{11i}}{u_{11}} - \frac{u^{ij}u_{11i} u_{11j}}{u_{11}^{2}}.
 \label{eq:Plog} \end{equation} 
 Now $$  u^{ij} u_{11ij} = \left( u^{ij} u_{11i}\right)_{j} - u^{ij}_{j}
 u_{11i}= -\left( u^{ij}_{1} u_{i1}\right)_{j}- u^{ij}_{j} u_{11i}. $$
 Expanding out the derivatives we get
 \begin{equation} u^{ij} u_{11ij} = - u^{ij}_{1j} u_{i1} + u^{ia} u^{jb} u_{1ab}u_{1ij}- u^{ij}_{j}
 u_{11i}. \label{eq:identity}\end{equation}
 Here we recognise the expression $-u^{ij}_{1j}u_{i1}$ as the co-efficient
  $G_{11}$ of the tensor $G$. (In fact the manipulation above is essentially
  the familiar identity, in Riemannian geometry,
   for $ \Delta \vert d\phi \vert^{2} $ involving the Ricci tensor,
   where $\phi$ is a harmonic function, but we have preferred to do the calculation
   directly.)
   Substituting (52) into (51), two terms cancel and we get
   \begin{equation} P \log(u_{11}) = \frac{G_{11}}{u_{11}} + \frac{1}{u_{11}} u^{ia}
   u^{jb} u_{1ab} u_{1ij} - \frac{1}{u_{11}^{2}} u^{ij} u_{11i} u_{11j}.
   \label{Plog2}\end{equation}
   Now simple calculations give
   \begin{eqnarray}   P(\frac{u_{1}^{2}}{2}) &=& u_{11},\\ 
   P (\log (-u)) &=& \frac{n-v^{i}u_{i}}{u} - \frac{u^{ij}u_{i} u_{j}}{u^{2}}.
   \label{Pids}\end{eqnarray}
     So we conclude that, at the maximum
   point,
   $$ \frac{ n- v^{i}u_{i}}{u} + u_{11} + \frac{G_{11}}{u_{11}} + (A) - (B) - (C)\leq
   0, $$
   where (A), (B) and (C) are the positive quantities
   \begin{eqnarray} (A)&=& \frac{1}{u_{11}} u^{ia} u^{jb} u_{1ab} u_{1ij}, \\
    (B) &=& \frac{1}{u_{11}^{2}} u^{ij} u_{11i} u_{iij}, \\
    (C)&=& \frac{1}{u^{2}} u^{ij} u_{i} u_{j}. \end{eqnarray}
We now use the condition (50) on the first derivatives. We may
suppose that  $u_{ij}$ is diagonal at the given maximum point, with
diagonal entries $u_{ii}=\lambda_{i}$ say. Then we have for $i\neq 1$
$$   \frac{u_{i}}{u} = - \frac{u_{11i}}{u_{11}}. $$
The term (C) is
$$  (C) = \frac{1}{u^{2}}\sum \lambda_{i}^{-1} u_{i}^{2}. $$
We write this as $$(C) = \frac{u_{1}^{2}}{ u^{2}u_{11}^{2}}+ (C)',
$$
where
$$ (C)'= \frac{1}{u^{2}} \sum_{i>1} \lambda_{i}^{-1} u_{i}^{2}= \sum_{i>1}
\frac{1}{\lambda_{i} \lambda_{1}^{2}} u_{11i}^{2}. $$
On the other hand, computing at this point, the other sums become
\begin{eqnarray} (A) &=& \sum_{ij} \frac{1}{\lambda_{1} \lambda_{i} \lambda_{j}} u_{1ij}^{2},
\\
 (B)&=& \sum_{i} \frac{1}{\lambda_{i} \lambda_{1}^{2}} u_{11i}^{2}. \end{eqnarray}
Once sees from this that
$$  (A)-(B)-(C)'\geq 0$$
so we conclude that at the maximum point
\begin{equation}\frac{n-u_{i}v^{i}}{u} + u_{11} +\frac{G_{11}}{u_{11}} - \frac{u_{1}^{2}}{u^{2}
u_{11}} \leq 0. \label{eq:conclusion}\end{equation}
This should be compared with the more standard calculation, for solutions
of the Monge-Amp\`ere equation $\det(u_{ij})=1$, when $v$ and $G$ vanish, so 
two terms (61) are absent. 

Now, computing the norm in the diagonal basis,
  $$\vert G\vert^{2} = \sum \lambda_{i}^{2} G_{ii}^{2} \geq \left( \frac{C_{11}}{u_{11}}\right)^{2},
  $$ so
  $$  \vert \frac{G_{11}}{u_{11}} \vert \leq c_{1}. $$
  Clearly $ n-u_{i} v^{i} \leq  n + c_{2} c_{3}$. Then
   (61) gives
  $$ \frac{N}{u} + u_{11} - \frac{u_{1}^{2}}{u^{2} u_{11}} \leq 0, $$
  where $N= n +c_{2} c_{3} + c_{1} c_{0}$. Thus we can adapt the argument
  from the standard case, replacing $n$ by $N$. If $h$ is the function
  $$  h= \exp(\log (-u) + f) = \vert u \vert u _{11} \exp(\frac{u_{1}^{2}}{2}), $$
  we see that at the maximum point for $h$
  $$ h^{2} - N \exp(\frac{u_{1}^{2}}{2}) h - u_{1}^{2} \exp(u_{1}^{2})
  \leq 0,$$
 which gives 
 $$  h_{max} \leq \exp(\frac{u_{1}^{2}}{2}) \frac{1}{2} \left( N+
  \sqrt{N^{2}
 + 4 u_{1}^{2}}\right). $$
 This gives
 $$  u_{11} \leq \frac{K}{\vert u\vert}, $$
 where
 $$   K= \frac{1}{2} \exp(\frac{c_{3}^{2}}{2})
  \left( N + \sqrt{N^{2} + 4 c_{3}^{2}}\right). $$
  
  \subsection{Yang-Mills estimate}
  
  In this subsection we consider a solution of the equation  $S(u)=A$ with $A$ a 
  {\it constant}. Roughly speaking, we show that when the dimension $n$
  is $2$ and given uniform bounds on the Hessian of
  $u$, the $L^{2}$ norm of the tensor $F$ (the \lq\lq Yang-Mills functional'')
  controls the $L^{\infty}$ norm.
  
  \begin{prop}
  Suppose $u$ is a convex function on the unit disc $D$ in $\bR^{2}$, smooth
  up to the boundary,
   which satisfies
  the equation $S(u)=A$ where $A$ is constant. Set
  $$ {\cal E}= \int_{D} \vert F \vert^{2}. $$
   If the Hessian of $u$ is
  bounded above and below
  $$    K^{-1}\leq (u_{ij})\leq K$$ then
  $$  \vert G(0)\vert^{2} \leq \kappa ( \cE + \cE^{3}) $$
  and for any $p>1$
  $$  \int_{\frac{1}{2} D}\vert F \vert^{p} \leq \kappa_{p} (\cE + \cE^{3})^{p}$$
  where $\kappa$ depends only on $K$ and $\sup_{D} \vert v \vert$ and $\kappa_{p}$
  depends on $K$ and $p$. 
  \end{prop}
  The proof makes use of the Sobolev inequalities. In dimension $2$ we
  have
      $$ \Vert \phi \Vert_{L^{p}} \leq C_{p} \Vert \nabla \phi \Vert_{L^{2}}$$
      for any $p$ and compactly supported $\phi$ on $D$ (here all norms
      are the standard Euclidean ones). The proof can be modified to give
      a similar result in dimension $3$ and extended to give information
      when $n=4$ provided ${\cal E}$ is sufficiently small, in the manner
      of Uhlenbeck\cite{kn:Uhl} and, still more, Anderson \cite{kn:An}
      and Tian and 
      Viaclovsky\cite{kn:TV}. In the proof we make more
   use of the Riemannian metric $g$ on $D \times \bR^{2}$
  defined by the convex function $u$. 
  
  The idea of the proof is to exploit the fact that the tensors $F$ and $G$ satisfy quasi-linear
  elliptic equations. To derive these we use the interpretation  of
  these tensors as (essentially) the Riemannian curvature and Ricci tensor
  of the metric $g$ see also \cite{kn:TV}. (Of course there is no particular difficulty in deriving
  these equations directly, without explicit reference to 
  Riemannian and Kahler geometry, but the derivation involves manipulating sixth order derivatives
  of the function $u$.)  
  \begin{lem} Suppose that $g$ is a Kahler metric of constant scalar curvature.
  Then the Riemann tensor $\Riem$ and Ricci tensor $\Ric$ of $g$ satisfy
  $$  \nabla^{*} \nabla \Ric = \Riem * \Ric, $$
  $$ \nabla^{*} \nabla \Riem = \Riem * \Riem + 2 \nabla' \nabla'' \Ric, $$
  where $*$ denotes appropriate natural algebraic bilinear forms and
  $\nabla', \nabla''$ are the $(1,0)$ and $(0,1)$ components of the covariant
  derivative. 
  \end{lem}
  To derive these identities we can consider more generally the curvature
  tensor $\Phi $ of a holomorphic
  vector bundle $E$ over a Kahler manifold $M$. Then we have $\db$-operators
  $$ \db: \Omega^{p,q}(\End E) \rightarrow \Omega^{p,q+1}(\End E). $$
  The Laplacians $ \nabla^{*}\nabla$ and $\Delta_{\db}= 2(\db^{*} \db + \db \db^{*})$
  differ by a Weitzenbock formula involving the curvature of the bundle
  and base manifold. Since $\db \Phi=0$ we have
  $$ \Delta_{\db} \Phi = 2 \db \db^{*} \Phi. $$
  Now the Kahler identities give
  $$ \db^{*} \Phi = = i \nabla' ( \Lambda \Phi), $$
  where $\Lambda: \Omega^{1,1}(\End E) \rightarrow \Omega^{0}(\End E)$ is the trace on the
  form component. Thus $\nabla^{*}\nabla \Phi$ is equal to $2i \db \nabla(\Lambda
  \Phi)= 2i \nabla'' \nabla' (\Lambda \Phi)$ plus a bilinear algebraic
  term involving $\Phi$ and the curvature of the base manifold. 
  
  We apply this first to the bundle $E=\Lambda^{n} TM$, the anticanonical
  line bundle. In this case $\Phi$ is essentially the Ricci tensor and
  $\Lambda \Phi$ is the scalar curvature. So the constant scalar curvature
  condition gives $\Delta_{\db} \Phi=0$ and the first formula follows from
  the discussion above. The second formula does not use the constant scalar
  curvature condition. It follows from the discussion above taking $E=TM$,
  when $\Phi$ is the Riemann curvature tensor and $\Lambda \Phi$ is the
  Ricci tensor.
  
  In our situation we deduce that $F,G$ satisfy equations which we write,
  rather schematically as
  \begin{eqnarray} \nabla^{*} \nabla G &=& F* G\\
                   \nabla^{*} \nabla F &=& F*F + 2 \nabla' \nabla'' G \label{FGids}
                   \end{eqnarray}

  \  
  The strategy now is to take the $L^{2}$ inner product of these identities
   with suitable compactly-supported tensor fields. It is convenient
  here to regard $g$ as a metric on the $4$-manifold $D\times T^{2}$ where
  $T^{2}$ is the torus $\bR^{2}/ {\bf Z}^{2}$. However all our data will
  depend only on the $D$ variables and the torus factor plays an entirely
  passive role: the push-forward of the volume form on the $4$-manifold
  to $D$ is the standard Lebesgue meausure and the Laplacian of $g$, applied
  to $T$-invariant functions is the operator $P$. Notice the crucial fact
  that for fixed $K$ the metric
  $g$ is uniformly equivalent to the standard Euclidean metric. 
  
  To prove the first part  of Proposition 6 we fix a standard cut-off function $\beta$
   on the disc, equal to $1$ in a neighbourhood
  of $0$, and take the $L^{2}$ inner product with $\beta^{2} G$ on either
  side of (62). This yields
  $$   \int_{D\times T^{2}} \nabla(\beta^{2}G).\nabla G \leq C \int_{D} \beta^{2}
  \vert G \vert^{2}  \vert F\vert. $$
  (In this proof $C$ will denote an unspecified constant, changing from
  line to line.)
  We have $$ \nabla(\beta^{2} G) .\nabla G = \vert \nabla(\beta G) \vert^{2}-
  \vert G \vert^{2} \vert \nabla \beta \vert^{2}$$ so we get
  \begin{equation} \int_{D\times T^{2}} \vert \nabla (\beta G) \vert^{2} \leq C \int_{D} \vert
  G \vert^{2} \sqrt \nabla \beta \vert^{2} + \beta^{2} \vert G \vert^{2} \vert
  F \vert. \label{gradbound}\end{equation}
  Applying the Sobolev inequality to $\vert \beta G \vert$ and Cauchy-Schwartz
   we obtain
  $$ \Vert \beta G \Vert^{2}_{L^{6}} \leq C \left( \Vert G \Vert_{L^{2}}^{2}
  + \Vert \beta G \Vert_{L^{4}}^{2} \Vert f \Vert_{L^{2}}\right). $$
  (Here we are using the fact that, for fixed $K$,
   the metric $g$ is uniformly equivalent to the Euclidean metric, so we
   can transfer the standard Sobolev inequalities to our setting.)
  Using 
  $$ \Vert \beta G \Vert_{L^{4}}^{4} \leq \Vert \beta G \Vert_{L^{2}} \Vert
  \beta G \Vert_{L^{6}}^{3}, $$
  we get
  $$ \Vert \beta G \Vert_{L^{6}}^{2} \leq C\left( \Vert G \Vert_{L^{2}}^{2}
  + \Vert \beta G \Vert_{L^{6}}^{3/2} \Vert F \Vert_{L^{2}} \Vert \beta
  G \Vert_{L^{2}}^{1/2}\right) \leq C\left( \cE^{2} +\cE^{3/2} \Vert \beta
  G \Vert_{L^{6}}^{3/2}\right). $$
  This yields
  $$ \Vert \beta G \Vert_{L^{6}} \leq C ( \cE + \cE^{3}). $$
  Thus we have an $L^{6}$ bound on $ G$ in a neighbourhood of $0$ and this
  gives an $L^{3/2}$ bound on $F* G$ over this neighbourhood. Going back to
   (64), we also have an $L^{2}$ bound on
  the derivative of $\beta G$:
  $$ \Vert \nabla(\beta G) \Vert_{L^{2}} \leq C ( \cE + \cE^{3}). $$
  We introduce another cut-off function $\gamma$, supported on the neighbourhood
  where $\beta =1$. Then
  we have
  $$ \Delta (\gamma \vert G \vert) \geq \vert \nabla^{*}\nabla (\gamma
  G) \vert \geq
  \vert \Delta \gamma \vert \vert G \vert + \vert \nabla \gamma \vert \vert \nabla
  G \vert + C \gamma \vert F \vert \vert G\vert. $$
  Now $$\Delta \gamma= u^{ij} \gamma_{ij} + v^{i} \gamma_{i}.$$
  So $\vert \Delta \gamma \vert $ satisfies a bound, depending on $\sup
  \vert v \vert$
  and  we have
  $$ \Delta (\gamma \vert G \vert) \geq \sigma, $$
  say where
  $$ \Vert \sigma \Vert_{L^{3/2}} \leq C (\cE+ \cE^{3}). $$
  Then we can apply Theorem 8.15 of \cite{kn:GT} (proved by the Moser iteration
  technique) to obtain the desired bound on
  $\vert G \vert$ at $0$. 
  
  For the second part of Proposition 6 we operate with the equation (62).
  We take the $L^{2}$ inner product with $\gamma^{2} F$ and integrate by
  parts. For the term involving the second derivatives of $G$ we write
  $$ \int_{D\times T^{2}} \gamma^{2} F . \nabla'' \nabla'G = \int_{D\times
  T^{2}} (\nabla'')^{*}
  (\gamma^{2} F ) .\nabla' G. $$
 This yields
 $$ \int_{D\times T^{2}} \vert \nabla(\gamma F) \vert^{2} \leq C \int_{D\times
 T^{2}} \gamma^{2} \vert F\vert^{3} + \vert F \vert^{2} \vert \nabla\gamma\vert^{2}
 + \gamma \vert \nabla(\gamma F) \vert \vert \nabla G + \gamma \vert \nabla
 \gamma \vert \vert F \vert \vert \nabla G \vert. $$

  Then, using the Sobolev inequlality as before and re-arranging we get
  an $L^{2}$ bound on the derivative of $\gamma F$ near $0$ which, in dimension
  $2$, gives the required $L^{p}$ bound (since we can suppose that $\gamma=1$
  on the disc $\frac{1}{2} D$).

  \subsection{Rescaling sections} 
  
  In this subsection we will bring together the three ingredients established
  above to obtain a pointwise bound on the tensor $G$ over compact subsets
  of $\Omega$.
  This involves rescaling the geometric data. A rescaling argument of this
  kind, using balls determined by the Riemannian distance function,
   would be fairly  standard. However there are difficulties in carrying
   this through unless one can establish some control of the injectivity
   radius, or something similar. We get around this difficulty by using
    Caffarelli's theory of the \lq\lq sections'' of a convex function, these
    taking the place of geodesic balls. 
    
    Recall from (5.1) that if $ u$ is a 
      smooth, strictly convex function  on an open set $U\subset
    \bR^{n}$ and  $x,y$ is a point  in $U$ we have defined
    $H_{x}(y)\geq 0$, vanishing if and only if $x=y$.
    For $t\geq 0$ the  {\it section} $S_{x}(t)$ at $x$ and  level $t$ is the
        set
        $$   S_{x}(t)= \{y\in U: H_{x}(y)\leq t\}. $$
        We will use three results about these sections, or equivalently
        the functions $H_{x}$, taken from \cite{kn:Gut}. For each of these results we suppose
        that the determinant of the
         Hessian satisfies upper and lower bounds
         $$ 0< \lambda \leq \det (u_{ij}) \leq \Lambda, $$
         throughout $U$, and the  constants $c_{i}$ below  depend only on
         $\lambda,\Lambda$. Recall that a  convex set $K$ in $\bR^{n}$
         is {\it normalised} if 
            \begin{equation}  \alpha_{n}B_{n} \subset  K \subset  B_{n} \label{eq:John}
            \end{equation}
            where $\alpha_{n}= n^{-3/2}$. Any compact convex
            set can be mapped to a normalised set by an affine-linear
             transformation (\cite{kn:Gut}, Theorem 1.8.2).

         \begin{prop}
         \begin{enumerate}
           \item {\rm (\cite{kn:Gut} page 50, Corollary 3.2.4)}  Suppose that $S_{x}(t)$ is compact. Then
            \begin{equation}
               c_{1} t^{n/2} \leq {\rm Vol}\ (S_{x}(t)) \leq c_{2} t^{n/2}.
             \end{equation} 
            \item {\rm (\cite{kn:Gut} page  55, Corollary 3.3.6 (i))} Suppose that $S_{x}(t)$ is compact and normalised. Then
            for $y\in S_{x}(t/2)$
            \begin{equation} d(y, \partial S_{x}(t)) \geq c_{3}>0. 
            \end{equation}
           
            \item  {\rm (\cite{kn:Gut} page 55, Theorem 3.3.7)} Suppose that
             $S_{x}(2t)$ is compact. Then if $H_{x}(y)\leq
            t$ and $H_{x}(z)\leq t$ we have $
            H_{y}(z)\leq c_{4} t $.             \item {\rm (\cite{kn:Gut} page 57, Theorem 3.3.8)} Suppose that
             $S_{x}(t)$ is compact and normalised. Then $S_{x}(t/2)$ contains
             the Euclidean ball of radius $c_{5}$ centred on $x$.
      \end{enumerate}\end{prop}
  
          The third result can be seen as a substitute for the triangle
          inequality, if one views $H_{x}(y)$ as a defining a notion of
          \lq\lq distance'' in $U$. (The assumption that $S_{x}(2t)$ be
          compact does not appear explicitly in \cite{kn:Gut},  where it is assumed
          that $U=\bR^{n}$ and all sections are compact, but a review of
          the proof shows that this is the hypothesis needed for our situation.)

         We will now discuss the scaling behaviour of the tensors we have
         associated to a convex function $u$. For $t>0$ and $T= T^{a}_{j}\in SL(n,
         \bR)$ we set
            $$ \tilde{u}(x) = t^{-1} u (\sqrt{t} x)\ ,\ u^{*}(x) = \tilde{u}(Tx)=t^{-1}u(\sqrt{t}
            T^{a}_{j}x^{j})
            $$
            We write $S= S^{j}_{a}$ for the inverse matrix of $T=T^{a}_{j}$.
            Then we have
            \begin{eqnarray}
             \det(u_{ij}) &=& \det(\tilde{u}_{ij}) = \det( u^{*}_{ij});\\
              \vert F^{*}\vert &=&  \vert \tilde{F}\vert = t \vert F
            \vert \\
             \vert G^{*}\vert&=& \vert \tilde{G}\vert = t \vert G\vert
  \\  \left(v^{*}\right)^{j} &=& S^{j}_{a} \tilde{v}^{a} =
            \sqrt{t} S^{j}_{a} v^{a}. 
            \label{eq:transform}\end{eqnarray}
             Here we have an obvious notation, in which for example
              $\tilde{v}$ and $ v^{*}$ refer to the vector fields obtained
              from the convex functions $\tilde{u}, u^{*}$. The verification
              of all these identities is completely elementary. 
              
              With these prelminaries in place we can proceed to our main
              argument. From now on we fix $n=2$ and suppose that
              $u$ is a convex function on $\Omega\subset \bR^{2}$ which
              satisfies Abreu's equation $S(u)=A$  with constant $A$. We
              set
              $$   E=\int_{\Omega} \vert F\vert^{2} $$  
              $$ \rho = \sup_{\Omega} \vert v\vert_{Euc}. $$
              Let $K\subset \subset K^{+}$ be compact subsets of $\Omega$
              and suppose that
              $$ 0< \lambda \leq \det(u_{ij}) \leq \Lambda, $$
              on $K^{+}$. For $x$ in $K$ we recall that
              $$  H_{x}(\partial K^{+}) =
               \min_{ y\in \partial K^{+}} H_{x}(y), $$
             and we let
             $$ \delta= H(K,K^{+})= \min_{x\in K} H_{x}(\partial K^{+}). $$
             We also put 
                $$ D= \max_{x\in K} H_{x}(\partial K). $$

                     Finally, we define a function $\Phi$ on $K$ by
                     $$ \Phi(x) = \vert G(x) \vert H_{x}(\partial K), $$
                     and let
                     $$  M= \max_{x \in K} \Phi(x). $$
                        
\begin{thm}
      There is a constant $\mu$ depending only on $\lambda,\Lambda, E,
      \rho, \delta, D$ such that $M\leq \mu$. \end{thm}

      To prove this  we may obviously suppose that $M>2$.
      We  consider a point $x_{0}$ where $\Phi$ attains its maximum
      value $M$ and set $t=\vert G(x_{0})\vert^{-1}$. Then $M> 2$ implies
      that $2t< H_{x_{0}}(\partial K)$ so the section 
      $\Sigma =S_{x_{0}}(t)$ lies in $K$. Hence $\Sigma$ is compact in
      $\Omega$. Our first goal is to show that $\vert G \vert$ over $\Sigma$
      is controlled by $t^{-1}$, i.e. by its value at $x_{0}$

      Now we may suppose that $$M\geq \max ( 2c_{4}^{2}, 4 c_{4} D \delta^{-1}),
      $$
      where $c_{4}>1$ is the constant of (3) in Proposition 7, depending on the given bounds
      $\lambda,\Lambda$.
      This means that if we define $\epsilon$ to be 
      $$ \epsilon= \min( 1/2,  \frac{c_{4} \delta}{2D}) $$
       we have $M\geq c_{4}^{2}/\epsilon$.

       We will now make two applications of the inequality in (3) of Proposition
       7.
       For the first we simply observe that in fact
        $S_{x_{0}}(2t)$ lies in $K$ and {\it a fortiori} in $K^{+}$. Hence
        this section is compact. Then for $y\in S_{x_{0}}(t)$
        we have $H_{x_{0}}(y)\leq t$ and trivially $H_{x_{0}}(x_{0})=0 \leq
        t$, so we deduce that $H_{y}(x_{0})\leq c_{4} t$.

       Now set $r=H_{x_{0}}(\partial K)$ so $t=r/M$. We claim that if
       $y$ is in $S_{x_{0}}(t)$ we have 
       \begin{equation} H_{y}(\partial K)\geq  \frac{r \epsilon}{c_{4}}. 
       \label{eq:crux}\end{equation}
       To see this, the result above yields
       $$ H_{y}(x_{0}) \leq c_{4}t = c_{4} r/M, $$
       and $c_{4}r/M\leq r \epsilon/c_{4}$ since we have arranged that $M \geq
       c_{4}^{2}/\epsilon$. Set $\tau= r\epsilon/c_{4}$. For $y\in S_{x_{0}}(t)\subset K$
       we have 
       \begin{equation}H_{y}(\partial K^{+}) \geq \delta \geq 
        \frac{2\epsilon D}{c_{4}}
       \geq \frac{2\epsilon r}{c_{4}}=2\tau.\label{eq:crux2} \end{equation}
       Suppose that (72) is not true, so there is a point $z\in \partial K$
       with $H_{y}(z) \leq \tau$. Then, by (73), $S_{y}(2\tau)$ lies in $K^{+}$
       and $H_{y}(x_{0}) \leq \tau$ so (3) in Proposition 7  would give
       $$  H_{x_{0}}(z) \leq c_{4} \tau = r\epsilon \leq r/2, $$
       a contradiction to $H_{x_{0}}(\partial K) =r$. 
       
       Now from $H_{y}(\partial K) \leq r\epsilon/c_{4}$ and the definition
       of $M$ we obtain our first goal:
       \begin{equation}   \vert G(y)\vert \leq 
        \frac{c_{4}}{\alpha} \vert G(x_{0})\vert 
       \label{eq:Gbound}\end{equation}
       for all $y\in \Sigma = S_{x_{0}}(t)$.

       \
       
       We now invoke our scaling construction for the restriction of $u$
       to $\Sigma$. We fix the real parameter $t$ to be as
       above. We know that there is some $k>0$ and an unimodular affine
       transformation $T$ so that $kt^{-1/2} T^{-1}(\Sigma)$ is normalised. But
       we know from (1) of Proposition 7 that the volume of $\Sigma$ lies between $c_{1} t
       $ and $c_{2}t$, so the volume of $kt^{-1/2} T^{-1}(\Sigma)$ lies between
       $k^{2} c_{1} $ and $k^{2} c_{2}$ hence
       \begin{equation}    \frac{\pi}{8c_{2}} \leq
       k^{2} \leq \frac{\pi}{c_{1}}. \label{eq:kbounds}\end{equation}
       Thus the convex set $\Sigma^{*}= t^{-1/2} T^{-1}(\Sigma)$ differs
       from its normalisation by a scale factor which is bounded above
       and below, so we can apply the results of (2) and (4) in Proposition
       7, with a change in the
       constants depending on the above bounds for $k$ (alternatively,
       we could normalise $\Sigma^{*}$ by changing the definition of $t$).
       \begin{lem} There is a constant $C$ depending only on $\lambda,
       \Lambda$ such that
       $$ \vert v^{*}\vert_{Euc} \leq C \vert v \vert_{Euc} {\rm Diam}(\Sigma)
       . $$
       \end{lem}
               
       To see this we may suppose that $T$ is diagonal, with $T^{1}_{1}=
       \lambda, T^{2}_{2}=\lambda^{-1}, T^{2}_{1}= T^{1}_{2}=0$. Then,
       by (71), 
       $v^{*}$ has components 
       $$ \left(v^{*}\right)^{1} = \sqrt{t} \lambda^{-1} v^{1}\ , \ \left(v^{*}\right)^{2}=
       \sqrt{t} \lambda v^{2}. $$
       But $\Sigma^{*}$ contains the disc of radius of radius 
       $R=k^{-1} 2^{-3/2}$ so $\Sigma=\sqrt{t} T(\Sigma^{*})$ contains
       an ellipse with semi-axes $\sqrt{t} \lambda R, \sqrt{t} R \lambda^{-1}$.
       In particular this ellipse is contained in $\Omega$ so
       $$ \sqrt{t} \lambda, \sqrt{t} \lambda^{-1} \leq \frac{{\rm Diam}(\Omega)}{2R}.
       $$
       Thus \begin{equation}
       \vert v^{*} \vert_{{\rm Euc}}\leq \frac{{\rm Diam}(\Omega)}{2R} 
       \vert v \vert_{{\rm Euc}}. \end{equation}

       \

       The scaling behaviour (70) of the tensor $G$ and the result (72)
         established
       above means that $\vert G^{*}\vert\leq c_{4}/\alpha$ over $\Sigma^{*}$.
       We now apply Proposition 5  to an  interior set
       $ \Sigma^{*}_{0}\subset \Sigma^{*}$. There is no loss of generality in supposing that
       $x_{0}=0$ and that $u$ vanishes to first order at this point, so $\Sigma^{*}= \{ y: u^{*}(y)\leq
       1\}$. We define 
       $$\Sigma^{*}_{0} = \{ y : u^{*}(y)\leq 3/4\}, $$
       and we set $u^{*}_{0}= u^{*}-3/4$. The lower bound on the distance to $\partial
       \Sigma^{*}$ furnished by (2) of Corollary 7 gives, in an elementary way,
        a bound on the derivative of $u^{*}_{0}$ over $\Sigma^{*}_{0}$
       and we can apply Proposition 5 to $u^{*}_{0}$, 
       taking $D=\Sigma^{*}_{0}$.
      Since we have bounds on $\vert G^{*}\vert$ and $\vert v^{*}\vert_{{\rm
      Euc}}$ (the latter by Lemma 8), we get an upper bound on the Hessian of $u^{*}_{0}$ over
   the further interior set
      $$\Sigma^{*}_{-1} = \{ y : u^{*}(y) \leq 1/2\}. $$
   
    The bound on $\det (u^{*}_{ij})$
      then  yields a lower bound on the Hessian over $\Sigma^{*}_{-1}$. We
      then use (4) of Proposition 7 to see that $\Sigma^{*}_{-1}$ contains a disc
       $\Delta$ about $0$ of fixed
      radius, and we can apply Proposition 6  to conclude that
      $$ \vert G^{*}(0)\vert^{2} \leq C (\cE + \cE^{3}), $$
      where $$\cE=\int_{\Delta} \vert F^{*}\vert^{2}
      \leq C \int_{\Sigma^{*}} \vert F^{*}\vert^{2}. $$
      (Here we use again the bound on $\vert v^{*}\vert_{{\rm Euc}}$.)
      Now the scaling behaviour of $F$ implies that
      $$ \int_{\Sigma^{*}} \vert F^{*} \vert^{2} = t \int_{\Sigma}\vert F\vert^{2}
      \leq t E. $$
      So we conclude that 
      $$ \vert G^{*}(0) \vert\leq C (\sqrt{tE + t^{3} E^{3}}, $$
      for some constant $C$ depending only on $\lambda,\Lambda, \delta,
      D, \rho$. But by construction $\vert G^{*}(0)\vert=1$, so we obtain
      a lower bound on $t$, thus an upper bound on $\vert G(x_{0}) \vert$ and
      thence on $M$. This completes the proof of Theorem 7.

      Given Theorem 7 it is straightforward to complete the proof of Theorem
      1
      in the case when $A$ is a constant and $\Omega\subset \bR^{2}$ is
      a polygon. By the discussion of uniform convexity in (5.1) above we
      can, for any compact set $K_{0}\subset \Omega$, find further compact
      sets
      $$  K_{0} \subset \subset K \subset \subset K^{+} $$
      in $\Omega$ such that
      $ \delta= H(K, K^{+})$ and $\delta_{0}=H(K_{0}, K)$ are bounded
below by positive bounds depending continuously on the data. Then Theorem
7 implies
that the tensor $G$ is bounded on  $K$. Covering $K$
by a finite number of sections and applying the same argument as above
we get upper and lower  bounds on the Hessian $(u_{ij})$ over a neighbourhood
of $K$. Further, we can apply the second part of Proposition 6 to get bounds on the $L^{p}$ norm of the
tensor $F$ for any $p$. Since $u_{ij}$ is bounded above and below these
are equivalent to $L^{p}$ bounds on the second derivatives of the matrix
$(u^{ij})$. In dimension $2$,  $L^{p}_{2}$ functions are continuous, so
 we deduce $L^{p}_{2}$ bounds on the inverse matrix $(u_{ij})$.
Thus we have uniform $L^{p}_{4}$---hence $C^{3,\alpha}$---bounds on
 $u$ over a neighbourhood of
$K$. From this point elementary methods suffice to bound all higher derivatives.

\end{document}